\pdfoutput=1
\documentclass[11pt]{article}
\usepackage[T1]{fontenc}
\usepackage{lmodern}
\usepackage[margin=1in]{geometry}
\usepackage{amsmath,amssymb,amsthm,mathtools}
\usepackage{microtype}
\usepackage{xurl}
\usepackage[authoryear,round]{natbib}
\setcitestyle{aysep={},yysep={;},notesep={, }}
\usepackage[hidelinks]{hyperref}
\hypersetup{
  pdftitle={Singleton Higson Coronas of Uncountable Groups},
  pdfauthor={Hussain Al-Rasheed},
  pdfkeywords={Higson corona, coarse geometry, periodic groups, locally nilpotent groups, locally finite groups, almost malnormal subgroups}
}
\newtheorem{theorem}{Theorem}[section]
\newtheorem{lemma}[theorem]{Lemma}
\newtheorem{proposition}[theorem]{Proposition}
\newtheorem{corollary}[theorem]{Corollary}
\theoremstyle{definition}
\newtheorem{definition}[theorem]{Definition}
\newtheorem{example}[theorem]{Example}
\theoremstyle{remark}
\newtheorem{remark}[theorem]{Remark}
\newtheorem{question}[theorem]{Question}

\newcommand{\FundingDisclosure}{No funding was received for this work.}
\newcommand{\InterestDisclosure}{The author declares that there are no competing interests.}

\title{Singleton Higson Coronas of Uncountable Groups}
\author{Hussain Al-Rasheed\thanks{General Studies Department, Jubail Industrial College, Royal Commission for Jubail and Yanbu, Jubail Industrial City 31961, Saudi Arabia. Corresponding author: \href{mailto:rashedhs@rcjy.edu.sa}{\nolinkurl{rashedhs@rcjy.edu.sa}}}}
\date{}

\begin{document}
\maketitle

\begin{abstract}
Equip a discrete group with the one-sided finitary coarse structure. Its Higson corona is a singleton precisely when every bounded slowly oscillating function is constant modulo a function vanishing at infinity. We give algebraic and geometric criteria for this property in uncountable groups. The main geometric result extends Holt's deleted coset-graph connectivity argument from periodic groups to nilpotent groups. It yields a singleton-corona theorem for every uncountable nonperiodic group whose two-generated subgroups are periodic or nilpotent. At cardinality $\aleph_1$, we show that a periodic group with nonsingleton corona must admit a continuous chain of countably infinite almost malnormal subgroups satisfying an additional restriction on commuting conjugates. This restriction proves that certain transfinite iterated restricted wreath products have singleton coronas, including explicit locally finite $p$-groups not covered by the preceding sufficient criteria. We also establish transporter and permanence criteria and reduce any counterexample among uncountable locally nilpotent groups to a locally finite $p$-group of cardinality $\aleph_1$. The proofs use countable subgroups containing translation defects and connectivity of coset graphs. We distinguish the new results from known cases, including groups with infinite center and periodic groups of cardinality greater than $\aleph_1$.
\end{abstract}

\par

\section{Introduction}\label{sec:intro}

Let $G$ be an infinite discrete group equipped with the right-translation finitary coarse structure, whose entourages are subsets of
\[
  E_K \coloneqq \{(x,y)\in G\times G : x^{-1}y\in K\},\qquad K\subseteq G\text{ finite}.
\]
The Higson corona $\partial G$ is the compact Hausdorff spectrum of $C_h(G)/C_0(G)$, where $C_h(G)$ consists of the bounded complex-valued functions whose right-translation defects $x\mapsto f(xg)-f(x)$ vanish at infinity for every $g\in G$; see \citep{Roe}. Thus $\partial G$ is a singleton exactly when every Higson function is constant modulo $C_0(G)$. This is stronger than connectedness of $\partial G$, which is equivalent to every right almost invariant subset of $G$ being finite or cofinite, or to $G$ having one coarse end.

For the balanced finitary coarse structure, which controls multiplication on both sides, the corona of every uncountable group is a singleton \citep[Theorem~1]{Protasov2019}. The one-sided structure depends on the group: uncountable free groups have disconnected coronas, whereas Holt's theorem gives connectedness for uncountable locally finite groups \citep{Holt1981,Cornulier,Protasov2005}. More generally, every uncountable periodic group is one-ended, by Holt's theorem and the splitting theorem for groups that are not locally finite \citep[Theorems~1.2 and~1.4]{Cornulier}. The question in this setting is therefore whether the connected full corona has more than one point.\footnote{Normality also distinguishes the countable and uncountable cases: the finitary coarse structure of a group is normal exactly when the group is countable. This is proved in T.~Banakh and I.~Protasov, \emph{The normality and bounded growth of balleans}, preprint (2018), \href{https://arxiv.org/abs/1810.07979}{arXiv:1810.07979}.} Outside the group setting, finitary coarse spaces can realize broad classes of prescribed Higson and binary coronas \citep{BanakhProtasov2020}.

Several singleton-corona results are already known. Protasov proved the uncountable abelian case \citep[Proposition~4]{Protasov2005}, and Filali and Protasov proved the case of uncountable $\omega_1$-normal groups, in which every countable subset lies in a countable normal subgroup \citep[Theorem~3.1]{FilaliProtasov2005}. Eastridge's vanishing theorem for $H^1(G,C_0(G))$ gives the same conclusion for uncountable groups with infinite center and for periodic groups of cardinality greater than $\aleph_1$ \citep[Theorem~7.9]{Eastridge2017}. Indeed, after inversion a Higson function $\alpha$ defines the cocycle $g\mapsto g\cdot\alpha-\alpha$, where $(g\cdot\alpha)(x)=\alpha(g^{-1}x)$. If this cocycle is $g\cdot b-b$ for some $b\in C_0(G)$, then $\alpha-b$ is invariant under all left translations and hence constant. These singleton-corona conclusions are therefore prior results. The criteria below address further cases of Protasov's detection question \citep[Question~1]{Protasov2019}.

The principal contributions are as follows.
\begin{enumerate}
  \item Theorem~\ref{thm:nilpotent_cosets} proves that, for a nilpotent group generated by two proper subgroups, deleting the two base vertices from the associated bipartite coset graph leaves a connected graph. This extends the periodic coset argument of Holt and Cornulier. The infinite dihedral example in Remark~\ref{rem:dihedral_sharp} shows that virtual abelianness does not suffice.
  \item Theorem~\ref{thm:nonperiodic_2gen} proves that every uncountable nonperiodic group whose two-generated subgroups are periodic or nilpotent has a singleton Higson corona. In particular, it applies to every uncountable nonperiodic locally nilpotent group, including those of cardinality $\aleph_1$.
  \item Lemma~\ref{lem:commuting_conjugates} gives an obstruction involving commuting conjugates. Together with almost-malnormality rigidity, it restricts the form of periodic counterexamples at $\aleph_1$ (Theorem~\ref{thm:aleph1_structure}). Theorem~\ref{thm:wreath} applies this obstruction to transfinite iterated restricted wreath products; Proposition~\ref{prop:wreath_separation} explains why the earlier sufficient criteria do not cover these groups.
\end{enumerate}

The method begins with a countable subgroup containing all translation defects associated with its own elements, called a \emph{defect hull}. The iterative subgroup construction occurs in Cornulier's account of Holt's method \citep[\S10]{Cornulier}; the support-closure construction and countable-support rigidity for general functions have antecedents in Eastridge's work \citep[Lemma~7.1 and Proposition~7.5]{Eastridge2017}. The author's exposition \citep{RashedHolt} presents the finite-defect construction and Holt's periodic coset argument in graph-theoretic form. We recall the shared ingredients and use them to develop the nilpotent connectivity theorem and the further rigidity criteria stated above.

Section~\ref{sec:prelim} gives the defect-hull construction and derives transporter rigidity. Section~\ref{sec:algebraic} gives the transporter and commuting-product criteria, applications to groups with countable supports, and permanence under finite index, suitable subgroup coverings, and uncountable ascendant subgroups. Finite-index invariance is standard; the proof is included for completeness.

Section~\ref{sec:aleph2} develops the coset-connectivity argument and gives a geometric proof of the singleton conclusion at cardinality at least $\aleph_2$ for groups whose two-generated subgroups are periodic or nilpotent. Its periodic case is already covered by Eastridge's theorem, and its nonperiodic case is strengthened in Section~\ref{sec:aleph1}. The independent geometric result is the nilpotent coset-connectivity theorem. Section~\ref{sec:aleph1} also proves the wreath-product application and reduces any locally nilpotent counterexample to a locally finite $p$-group of cardinality $\aleph_1$. Section~\ref{sec:questions} records the remaining questions; no converse to the counterexample restrictions is asserted.

The singleton condition also has useful analytic interpretations. The equality $C_h(G)=\mathbb C\cdot1+C_0(G)$ identifies the Higson compactification with the one-point compactification of the discrete space $G$. On $\ell^2(G)$, let $(M_f\xi)(x)=f(x)\xi(x)$ and $(\rho(g)\xi)(x)=\xi(xg)$. Since
\[
  ((\rho(g)M_f-M_f\rho(g))\xi)(x)
  =(f(xg)-f(x))\xi(xg),
\]
the commutators $[M_f,\rho(g)]$ are compact for all $g$ exactly when $f\in C_h(G)$. Hence $\partial G$ is a singleton exactly when every such diagonal operator is scalar \emph{modulo compact operators}:
\[
  [M_f,\rho(g)]\in\mathcal K(\ell^2(G))\quad(g\in G)
  \quad\Longrightarrow\quad
  M_f-cI\in\mathcal K(\ell^2(G))\quad\text{for some }c\in\mathbb C.
\]
Finally, a bounded function on $G$ extends to the Stone--\v{C}ech compactification $\beta G$, and its translation defects lie in $C_0(G)$ precisely when its extension is right $G$-invariant on $\beta G\setminus G$. Singletonness is therefore equivalent to constancy of every continuous right $G$-invariant function on that remainder.

This phenomenon does not occur for countably infinite groups. Such a group admits a proper left-invariant metric inducing its finitary coarse structure. Choose pairwise disjoint balls centered at points $x_n$, with radii tending to infinity. Any bounded sequence can be extended from the $x_n$ by Lipschitz bumps supported on these balls, with Lipschitz constants tending to zero. The resulting function is Higson, so restriction to the $x_n$ induces a surjection onto $\ell^\infty(\mathbb N)/C_0(\mathbb N)$. By Gelfand duality, the Higson corona contains a copy of $\beta\mathbb N\setminus\mathbb N$. Thus even one-ended countable groups have nonsingleton full coronas.

\medskip\noindent\textbf{Conventions.} All groups carry the discrete topology, and $e$ denotes the identity element. Throughout, \emph{countable} means at most countable (cardinality at most $\aleph_0$). For any set $X$, we write $C_0(X)$ for the space of functions $f\colon X\to\mathbb{C}$ such that $\{x\in X : |f(x)|\ge\varepsilon\}$ is finite for every $\varepsilon>0$. For subsets $A,B$ of a set $X$, we write $A\triangle B \coloneqq (A\setminus B)\cup(B\setminus A)$ and $\chi_A\colon X\to\{0,1\}$ for the characteristic function of $A$. For elements $x,y$ of a group $G$, we write $[x,y]\coloneqq xyx^{-1}y^{-1}$, and for subsets $S,H\subseteq G$ we define the \emph{transporter} $T_G(S,H)\coloneqq\{t\in G : t^{-1}St\subseteq H\}$.

\section{Preliminaries and countable defect hulls}\label{sec:prelim}

Let $G$ be an infinite group. A bounded function $f\colon G\to\mathbb{C}$ is \emph{slowly oscillating}, or \emph{Higson}, if
\begin{equation}\label{eq:translate}
  \bigl(x\longmapsto f(xg)-f(x)\bigr)\in C_0(G)\qquad\text{for every }g\in G.
\end{equation}
We write $C_h(G)$ for the unital $C^*$-subalgebra of bounded Higson functions on $G$. Since $G$ is infinite, $1\notin C_0(G)$, so $C_h(G)/C_0(G)$ is a nonzero commutative unital $C^*$-algebra. Its compact Hausdorff Gelfand spectrum is the \emph{Higson corona} $\partial G$, satisfying $C(\partial G)\cong C_h(G)/C_0(G)$. Consequently, $\partial G$ is a singleton if and only if every $f\in C_h(G)$ is constant modulo $C_0(G)$, meaning that $f-c\in C_0(G)$ for some $c\in\mathbb{C}$.

\begin{lemma}[Oscillation and cluster points]\label{lem:oscillation}
Let $X$ be an infinite set and let $f\colon X\to\mathbb{C}$ be a bounded function that is not constant modulo $C_0(X)$.
\begin{enumerate}
  \item There exist $\varepsilon>0$ and sequences $(u_n)_{n\ge 1}$ and $(v_n)_{n\ge 1}$ in $X$ such that all elements $u_1,v_1,u_2,v_2,\ldots$ are pairwise distinct and $|f(u_n)-f(v_n)|\ge\varepsilon$ for every $n\ge 1$.
  \item There exist distinct complex numbers $c_1\neq c_2$ and sequences $(u_n)_{n\ge 1}$ and $(v_n)_{n\ge 1}$ in $X$ whose combined terms $u_1,v_1,u_2,v_2,\ldots$ are pairwise distinct such that $f(u_n)\to c_1$ and $f(v_n)\to c_2$ as $n\to\infty$.
\end{enumerate}
\end{lemma}

\begin{proof}
For (1), we first claim that there exists $\varepsilon>0$ such that outside every finite subset of $X$ there exist $u,v$ with $|f(u)-f(v)|\ge\varepsilon$. Otherwise, for each integer $m\ge 1$ there exists a finite subset $F_m\subseteq X$ such that $|f(u)-f(v)|<1/m$ for all $u,v\in X\setminus F_m$. Choosing $x_m\in X\setminus(F_1\cup\cdots\cup F_m)$ gives a Cauchy sequence $(f(x_m))$ in $\mathbb{C}$ whose limit $c\in\mathbb{C}$ satisfies $|f(x)-c|\le 1/m$ for all $x\in X\setminus F_m$, whence $f-c\in C_0(X)$, a contradiction. Recursively choosing $(u_n,v_n)$ outside $\{u_1,v_1,\ldots,u_{n-1},v_{n-1}\}$ with $|f(u_n)-f(v_n)|\ge\varepsilon$ automatically gives $u_n\neq v_n$, so all terms $u_1,v_1,u_2,v_2,\ldots$ are pairwise distinct.

For (2), take $\varepsilon>0$ and $(u_n),(v_n)$ from (1). Since $f$ is bounded, the sequence $((f(u_n),f(v_n)))_{n\ge 1}$ is bounded in $\mathbb{C}^2$ and therefore admits a convergent subsequence $((f(u_{n_k}),f(v_{n_k})))\to(c_1,c_2)$. Because $|f(u_{n_k})-f(v_{n_k})|\ge\varepsilon$, we have $|c_1-c_2|\ge\varepsilon>0$, so $c_1\neq c_2$. Relabeling these subsequences yields (2).
\end{proof}

The following support observation is the bounded-function form of \citep[Proposition~7.5]{Eastridge2017}, with the translation convention adjusted by inversion.

\begin{lemma}[Countable-support rigidity]\label{lem:countable_support}
Let $G$ be an uncountable group, let $f\in C_h(G)$, and let $c\in\mathbb{C}$. If $\{x\in G : f(x)\neq c\}$ is countable, then $f-c\in C_0(G)$.
\end{lemma}

\begin{proof}
Set $E\coloneqq\{x\in G : f(x)\neq c\}$. Since $E$ is countable, $E^{-1}E$ is countable, so we may choose $h_0\in G\setminus E^{-1}E$. Then $Eh_0\cap E=\varnothing$, whence $f(xh_0)=c$ for every $x\in E$. For every $\varepsilon>0$,
\[
  \{x\in G : |f(x)-c|\ge\varepsilon\}
  = \{x\in E : |f(x)-f(xh_0)|\ge\varepsilon\}
  \subseteq \{x\in G : |f(xh_0)-f(x)|\ge\varepsilon\},
\]
which is finite by \eqref{eq:translate}. Hence $f-c\in C_0(G)$.
\end{proof}

As in the coset-graph proof of Holt's theorem, passing from right translations to left translations aligns slow oscillation with constancy on right cosets.

\begin{lemma}[Inversion and left-translation defects]\label{lem:inversion}
For a function $f\colon G\to\mathbb{C}$, define $\tilde{f}\colon G\to\mathbb{C}$ by $\tilde{f}(x)\coloneqq f(x^{-1})$. Then $f\in C_h(G)$ if and only if $\tilde{f}$ is bounded and, for every $g\in G$ and $\varepsilon>0$, the set
\[
  D_{g,\varepsilon}(\tilde{f})
  \coloneqq \{x\in G : |\tilde{f}(gx)-\tilde{f}(x)|\ge\varepsilon\}
\]
is finite. Moreover, $f$ is constant modulo $C_0(G)$ if and only if $\tilde{f}$ is constant modulo $C_0(G)$.
\end{lemma}

\begin{proof}
For every $g,x\in G$, we have $\tilde{f}(gx)-\tilde{f}(x)=f(x^{-1}g^{-1})-f(x^{-1})$. Since inversion is a bijection of $G$, it maps $D_{g,\varepsilon}(\tilde{f})$ onto $\{y\in G : |f(yg^{-1})-f(y)|\ge\varepsilon\}$, and for any $c\in\mathbb{C}$ it maps $\{x\in G : |\tilde{f}(x)-c|\ge\varepsilon\}$ onto $\{y\in G : |f(y)-c|\ge\varepsilon\}$.
\end{proof}

We call a bounded function $\alpha\colon G\to\mathbb{C}$ \emph{left-slowly oscillating} if $D_{g,\varepsilon}(\alpha)$ is finite for every $g\in G$ and $\varepsilon>0$. For such a function $\alpha$ and each $g\in G$, the \emph{exact defect set}
\[
  D_g(\alpha) \coloneqq \{x\in G : \alpha(gx)\neq\alpha(x)\}
  = \bigcup_{m=1}^{\infty} D_{g,1/m}(\alpha)
\]
is countable.

\begin{definition}\label{def:defect_hull}
Let $\alpha\colon G\to\mathbb{C}$ be left-slowly oscillating. A subgroup $M\le G$ is a \emph{defect hull} for $\alpha$ if $D_h(\alpha)\subseteq M$ for every $h\in M$.
\end{definition}

If $M\le G$ is a defect hull for $\alpha$, then for every $x\in G\setminus M$ and every $h\in M$ we have $x\notin D_h(\alpha)$, so $\alpha(hx)=\alpha(x)$. Thus $\alpha$ is constant on every right $M$-coset other than $M$.

The next lemma records the iterative closure construction just discussed. For characteristic functions, the exact defect sets are finite, as in \citet{RashedHolt}; here they are countable unions of finite level-defect sets, as in the support argument of \citet[Lemma~7.1]{Eastridge2017}. We include the proof to specify the least hull and the notation $\operatorname{cl}_\alpha$ used below.

\begin{lemma}[Canonical countable defect hulls]\label{lem:defect_hull}
Let $\alpha\colon G\to\mathbb{C}$ be left-slowly oscillating. For every subset $C\subseteq G$, there is a unique smallest defect hull $\operatorname{cl}_\alpha(C)\le G$ containing $C$. If $C$ is countable, then $\operatorname{cl}_\alpha(C)$ is countable.
\end{lemma}

\begin{proof}
Define an increasing chain of subgroups inductively by $M_0\coloneqq\langle C\rangle$ and
\[
  M_{n+1} \coloneqq \left\langle M_n\cup\bigcup_{h\in M_n} D_h(\alpha)\right\rangle \qquad (n\ge 0),
\]
and set $M\coloneqq\bigcup_{n=0}^{\infty} M_n$. Then $M$ is a subgroup containing $C$. If $h\in M$, then $h\in M_n$ for some $n$, so $D_h(\alpha)\subseteq M_{n+1}\subseteq M$; hence $M$ is a defect hull. Conversely, any defect hull $N\le G$ containing $C$ contains $M_0$ and, by induction, contains every $M_n$, hence $M\subseteq N$. Thus $M=\operatorname{cl}_\alpha(C)$ is the unique smallest defect hull containing $C$. Finally, if $C$ is countable, induction shows that each $M_n$ is countable (since each $D_h(\alpha)$ is countable), so $M$ is countable.
\end{proof}

\begin{lemma}[Defect-hull transporter rigidity]\label{lem:self_normalizing}
Let $G$ be a group, let $\alpha\colon G\to\mathbb{C}$ be left-slowly oscillating, and let $M<G$ be a defect hull for $\alpha$.
\begin{enumerate}
  \item For every $x\in G\setminus M$, the function $h\mapsto\alpha(x^{-1}hx)-\alpha(x)$ belongs to $C_0(M)$.
  \item If $M_0\subseteq G$ is a subset such that $\alpha|_{M_0}$ is not constant modulo $C_0(M_0)$, then $T_G(M_0,M)\subseteq M$. In particular, if $\alpha|_M$ is not constant modulo $C_0(M)$, then $N_G(M)=M$.
\end{enumerate}
\end{lemma}

\begin{proof}
For (1), fix $x\in G\setminus M$ and $\varepsilon>0$. Since $M$ is a defect hull, $\alpha(hx)=\alpha(x)$ for every $h\in M$, and therefore
\[
  \{h\in M : |\alpha(x^{-1}hx)-\alpha(x)|\ge\varepsilon\}
  = \{h\in M : |\alpha(x^{-1}(hx))-\alpha(hx)|\ge\varepsilon\}.
\]
Because $h\mapsto hx$ is injective on $M$ and $D_{x^{-1},\varepsilon}(\alpha)$ is finite, this set is finite, proving (1).

For (2), suppose for contradiction that there exists $t\in T_G(M_0,M)\setminus M$. Then $x\coloneqq t^{-1}\in G\setminus M$ satisfies $xM_0 x^{-1}=t^{-1}M_0 t\subseteq M$. As $m$ ranges over $M_0$, the conjugates $h\coloneqq xmx^{-1}$ range injectively over a subset of $M$, and $\alpha(x^{-1}hx)=\alpha(m)$. By (1), the function $m\mapsto\alpha(m)-\alpha(x)$ belongs to $C_0(M_0)$, contradicting the hypothesis on $M_0$. Taking $M_0=M$ gives $N_G(M)\subseteq T_G(M,M)\subseteq M$, so $N_G(M)=M$.
\end{proof}

\begin{corollary}\label{cor:self_normalizing}
Let $G$ be an uncountable group that contains no proper countably infinite self-normalizing subgroup (in particular, this holds if $G$ satisfies the normalizer condition, i.e., $N_G(H)>H$ for every proper subgroup $H<G$). Then $\partial G$ is a singleton.
\end{corollary}

\begin{proof}
Suppose for contradiction that some $f\in C_h(G)$ is not constant modulo $C_0(G)$, and set $\alpha(x)\coloneqq f(x^{-1})$. By Lemmas~\ref{lem:inversion} and~\ref{lem:oscillation}(1), there is a countably infinite subset $S\subseteq G$ on which $\alpha|_S$ is not constant modulo $C_0(S)$. Set $M\coloneqq\operatorname{cl}_\alpha(S)$. By Lemma~\ref{lem:defect_hull}, $M$ is a countably infinite defect hull (hence proper in $G$, since $G$ is uncountable) and $\alpha|_M$ is not constant modulo $C_0(M)$. Lemma~\ref{lem:self_normalizing}(2) then yields $N_G(M)=M$, contradicting the hypothesis.
\end{proof}

\section{Transporter, commuting-product, and permanence criteria}\label{sec:algebraic}

\begin{theorem}[Transporter criterion]\label{thm:transporter}
Let $G$ be a group. Suppose that for every countable subset $S\subseteq G$, there exists a countable subset $H\subseteq G$ such that $T_G(S,H)$ is uncountable. Then the Higson corona $\partial G$ is a singleton.
\end{theorem}

\begin{proof}
Taking $S=\{e\}$ shows that $G$ contains an uncountable subset and is therefore uncountable. Suppose for contradiction that some $f\in C_h(G)$ is not constant modulo $C_0(G)$, and set $\alpha(x)\coloneqq f(x^{-1})$. By Lemmas~\ref{lem:inversion} and~\ref{lem:oscillation}(1), there exists a countable subset $S\subseteq G$ such that $\alpha|_S$ is not constant modulo $C_0(S)$. By hypothesis, there is a countable subset $H\subseteq G$ with $T_G(S,H)$ uncountable; put $M\coloneqq\operatorname{cl}_\alpha(H)$. By Lemma~\ref{lem:defect_hull}, $M$ is a countable (hence proper) defect hull for $\alpha$ containing $H$, so Lemma~\ref{lem:self_normalizing}(2) yields $T_G(S,H)\subseteq T_G(S,M)\subseteq M$, contradicting that $T_G(S,H)$ is uncountable.
\end{proof}

For an infinite set $I$ and a field $F$, we write $\operatorname{SL}_{\mathrm{fin}}(I,F)$ for the group of $I\times I$ matrices over $F$ that differ from the identity matrix in only finitely many entries and have determinant $1$. For a totally ordered set $(Q,<)$ and a field $F$, the \emph{McLain group} $M(Q,F)$ is the group of $Q\times Q$ matrices over $F$ generated by the elementary transvections $I+\lambda e_{pq}$ ($p<q$ in $Q$, $\lambda\in F$).

\begin{corollary}\label{cor:transporter_apps}
Let $G$ be a group. In each of the following cases, $\partial G$ is a singleton:
\begin{enumerate}
  \item For every countable subset $S\subseteq G$, the centralizer $C_G(S)$ is uncountable. In particular, this holds for $\operatorname{Alt}_{\mathrm{fin}}(I)$, $\operatorname{Sym}_{\mathrm{fin}}(I)$, and $\operatorname{SL}_{\mathrm{fin}}(I,F)$ for any uncountable set $I$ and any field $F$, and for the McLain group $M(Q,F)$ over any uncountable totally ordered set $Q$ and any field $F$.
  \item $G$ is uncountable and $\omega_1$-normal, i.e., every countable subset of $G$ is contained in a countable normal subgroup of $G$ \textup{(\citep[Theorem~3.1]{FilaliProtasov2005})}.
\end{enumerate}
\end{corollary}

\begin{proof}
For (1), take $H\coloneqq S$; then $C_G(S)\subseteq T_G(S,S)$ is uncountable. In each listed example, every countable subset $S$ involves only countably many indices in $I$ (or $Q$), leaving an uncountable centralizer supported on the complement. For (2), let $H\trianglelefteq G$ be a countable normal subgroup containing $S$; then $T_G(S,H)=G$ is uncountable. Apply Theorem~\ref{thm:transporter}.
\end{proof}

A second criterion concerns commuting products and gives a direct proof of the known infinite-center case.

\begin{theorem}[Commuting-product criterion]\label{thm:commuting_product}
Let $G$ be a group containing subgroups $A,B\le G$ such that $G=AB$, $[A,B]=\{e\}$, $A$ is infinite, and $B$ is uncountable. Then $\partial G$ is a singleton. In particular, every uncountable group with an infinite center has a singleton Higson corona.
\end{theorem}

\begin{proof}
Fix $f\in C_h(G)$ and choose a countably infinite subset $A_0\subseteq A$. The union
\[
  E \coloneqq \bigcup_{a_0\in A_0,\, j\ge 1} \{x\in G : |f(xa_0)-f(x)|\ge 1/j\}
\]
is countable by \eqref{eq:translate}, and $f(xa_0)=f(x)$ for every $x\in G\setminus E$ and every $a_0\in A_0$.

Fix $x\in G$, and let $b_1,b_2\in B$ satisfy $xb_1,xb_2\notin E$. Put $h\coloneqq b_1^{-1}b_2\in B$. Because $[A,B]=\{e\}$, for every $a_0\in A_0$ we have $((xb_1)a_0)h=(xb_1 h)a_0=(xb_2)a_0$, and therefore
\[
  \bigl|f\bigl(((xb_1)a_0)h\bigr)-f\bigl((xb_1)a_0\bigr)\bigr|
  = |f((xb_2)a_0)-f((xb_1)a_0)|
  = |f(xb_2)-f(xb_1)|.
\]
Since $a_0\mapsto(xb_1)a_0$ is injective on the infinite set $A_0$, \eqref{eq:translate} applied to $h$ forces $f(xb_2)=f(xb_1)$.

For each $a\in A$, set $B_a\coloneqq\{b\in B : ab\notin E\}=B\setminus(a^{-1}E\cap B)$. Because $B$ is uncountable and $E$ is countable, $B_a$ is uncountable; by the preceding paragraph applied to $x=a$, the map $b\mapsto f(ab)$ takes a single constant value $c(a)\in\mathbb{C}$ on $B_a$. Moreover, $B_e\cap B_a$ is uncountable (hence infinite). For every $b\in B_e\cap B_a$, we have $ba=ab$, so
\[
  |f(ba)-f(b)| = |f(ab)-f(b)| = |c(a)-c(e)|.
\]
Applying \eqref{eq:translate} to $a$ gives $c(a)=c(e)$. Set $c\coloneqq c(e)$.

Given any $y\in G\setminus E$, write $y=ab$ with $a\in A$ and $b\in B$; then $b\in B_a$, so $f(y)=f(ab)=c(a)=c$. Thus $\{y\in G : f(y)\neq c\}\subseteq E$ is countable, and Lemma~\ref{lem:countable_support} yields $f-c\in C_0(G)$. Hence $\partial G$ is a singleton. The final assertion follows by taking $A\coloneqq Z(G)$ and $B\coloneqq G$.
\end{proof}

\begin{remark}\label{rem:eastridge_center}
The infinite-center conclusion is a special case of \citep[Theorem~7.9]{Eastridge2017}; the preceding proof treats the more general commuting-product hypothesis directly.
\end{remark}

Next we consider finite index, subgroup coverings with uncountable intersection, and uncountable ascendant subgroups. Finite-index invariance is a standard consequence of coarse equivalence; we include a direct proof in terms of Higson functions.

\begin{proposition}[Finite-index invariance]\label{prop:finite_index}
Let $G$ be an infinite group and $H\le G$ a subgroup of finite index. Then $\partial G$ is a singleton if and only if $\partial H$ is a singleton.
\end{proposition}

\begin{proof}
Since $G$ is infinite and $[G:H]<\infty$, the subgroup $H$ is infinite. Choose a finite set $R\subseteq G$ of right coset representatives with $e\in R$ so that $G=\bigsqcup_{r\in R} Hr$.

Suppose first that $\partial H$ is a singleton, and let $f\in C_h(G)$. For each $r\in R$, define $f_r\colon H\to\mathbb{C}$ by $f_r(h)\coloneqq f(hr)$. For any $k\in H$, we have
\[
  f_r(hk)-f_r(h) = f\bigl((hr)(r^{-1}kr)\bigr)-f(hr).
\]
Since $h\mapsto hr$ is injective on $H$ and $f\in C_h(G)$, it follows that $f_r\in C_h(H)$, so $f_r-c(r)\in C_0(H)$ for a unique $c(r)\in\mathbb{C}$. For any $r,s\in R$, writing $f_s(h)-f_r(h)=f\bigl((hr)(r^{-1}s)\bigr)-f(hr)$ shows that $f_s-f_r\in C_0(H)$, whence $c(s)-c(r)\in C_0(H)$. Since $H$ is infinite, $c(r)=c$ is independent of $r\in R$. Thus, for every $\varepsilon>0$,
\[
  \{x\in G : |f(x)-c|\ge\varepsilon\}
  = \bigcup_{r\in R} \{h\in H : |f_r(h)-c|\ge\varepsilon\}r
\]
is a finite union of finite sets, so $f-c\in C_0(G)$.

Conversely, suppose that $\partial G$ is a singleton, and let $f_0\in C_h(H)$. Extend $f_0$ to a bounded function $f\colon G\to\mathbb{C}$ by setting $f(hr)\coloneqq f_0(h)$ for $h\in H$ and $r\in R$. Fix $g\in G$. For each $r\in R$, write $rg=k_r s_r$ with $k_r\in H$ and $s_r\in R$. Then for every $h\in H$,
\[
  f((hr)g)-f(hr) = f(hk_r s_r)-f(hr) = f_0(hk_r)-f_0(h).
\]
Because $R$ is finite and $f_0\in C_h(H)$, we obtain $(x\mapsto f(xg)-f(x))\in C_0(G)$, so $f\in C_h(G)$. Hence $f-c\in C_0(G)$ for some $c\in\mathbb{C}$. Since $e\in R$, we have $f|_H=f_0$, and therefore $f_0-c\in C_0(H)$.
\end{proof}

\begin{proposition}[Union permanence]\label{prop:union}
Let $G$ be a group covered by a family of subgroups $\{N_i : i\in I\}$ such that $U\coloneqq\bigcap_{i\in I} N_i$ is uncountable. If $\partial N_i$ is a singleton for every $i\in I$, then $\partial G$ is a singleton.
\end{proposition}

\begin{proof}
Fix $f\in C_h(G)$. For each $i\in I$, the restriction $f|_{N_i}$ belongs to $C_h(N_i)$, so there exists $c_i\in\mathbb{C}$ with $f|_{N_i}-c_i\in C_0(N_i)$. Because $U\subseteq N_i$ is infinite, $f|_U-c_i\in C_0(U)$ for every $i\in I$, which forces $c_i=c$ to be independent of $i$.

Suppose for contradiction that $f-c\notin C_0(G)$. Then there exist $\varepsilon>0$ and a sequence $(z_k)_{k\ge 1}$ of pairwise distinct elements of $G$ with $|f(z_k)-c|\ge\varepsilon$ for all $k\ge 1$. Since $G=\bigcup_{i\in I} N_i$, choose for each $k\ge 1$ an index $i_k\in I$ with $z_k\in N_{i_k}$, and set
\[
  B_k \coloneqq \{u\in U : |f(z_k u)-c|\ge\varepsilon/2\}.
\]
Because $z_k\in N_{i_k}$ and $U\subseteq N_{i_k}$, left multiplication by $z_k$ maps $U$ injectively into $N_{i_k}$. Since $f|_{N_{i_k}}-c\in C_0(N_{i_k})$, each $B_k$ is finite. As $U$ is uncountable, we may pick $u_0\in U\setminus\bigcup_{k=1}^{\infty} B_k$, giving $|f(z_k u_0)-c|<\varepsilon/2$ for every $k\ge 1$. On the other hand, the distinct elements $z_k$ eventually leave every finite subset of $G$, so \eqref{eq:translate} applied to $u_0$ yields $|f(z_k u_0)-f(z_k)|<\varepsilon/2$ for all sufficiently large $k$. For such $k$,
\[
  |f(z_k)-c| \le |f(z_k)-f(z_k u_0)| + |f(z_k u_0)-c| < \frac{\varepsilon}{2}+\frac{\varepsilon}{2}=\varepsilon,
\]
a contradiction. Thus $f-c\in C_0(G)$.
\end{proof}

Recall that a subgroup $N\le G$ is \emph{ascendant} in $G$ if there exist an ordinal $\lambda$ and an ascending chain of subgroups $(N_\alpha)_{\alpha\le\lambda}$ with $N_0=N$, $N_\lambda=G$, $N_\alpha\trianglelefteq N_{\alpha+1}$ for all $\alpha<\lambda$, and $N_\beta=\bigcup_{\alpha<\beta} N_\alpha$ for every limit ordinal $\beta\le\lambda$.

\begin{theorem}[Ascendant propagation]\label{thm:ascendant}
Let $G$ be a group containing an uncountable ascendant subgroup $N$. If $\partial N$ is a singleton, then $\partial G$ is a singleton.
\end{theorem}

\begin{proof}
Let $(N_\alpha)_{\alpha\le\lambda}$ be an ascending series from $N_0=N$ to $N_\lambda=G$. We prove by transfinite induction on $\alpha\le\lambda$ that $\partial N_\alpha$ is a singleton. The base case $\alpha=0$ holds by hypothesis, and every limit step $\beta\le\lambda$ follows directly from Proposition~\ref{prop:union} applied to $N_\beta=\bigcup_{\alpha<\beta} N_\alpha$ (since $N_0\subseteq\bigcap_{\alpha<\beta} N_\alpha$ is uncountable).

It remains to verify the successor step: assuming $N\trianglelefteq G$ is uncountable with $\partial N$ a singleton, we prove that $\partial G$ is a singleton. Fix $f\in C_h(G)$. For each $x\in G$, define $f_x\colon N\to\mathbb{C}$ by $f_x(n)\coloneqq f(nx)$. As in the proof of Proposition~\ref{prop:finite_index}, for every $h\in N$ the identity
\[
  f_x(nh)-f_x(n) = f\bigl((nx)(x^{-1}hx)\bigr)-f(nx)
\]
shows that $f_x\in C_h(N)$, so there is a unique $c(x)\in\mathbb{C}$ with $f_x-c(x)\in C_0(N)$. For any $x,y\in G$, the identity $f_y(n)-f_x(n)=f\bigl((nx)(x^{-1}y)\bigr)-f(nx)$ gives $f_y-f_x\in C_0(N)$, and since $N$ is infinite, $c(x)=c(y)\eqqcolon c$.

Suppose for contradiction that $f-c\notin C_0(G)$, and choose $\varepsilon>0$ and a sequence $(z_k)_{k\ge 1}$ of pairwise distinct elements of $G$ with $|f(z_k)-c|\ge\varepsilon$. For each $k\ge 1$, set
\[
  B_k \coloneqq \{n\in N : |f(z_k n)-c|\ge\varepsilon/2\}.
\]
Because $N\trianglelefteq G$, conjugation $\varphi_k(n)\coloneqq z_k n z_k^{-1}$ is a bijection of $N$ satisfying $f(z_k n)=f(\varphi_k(n)z_k)=f_{z_k}(\varphi_k(n))$. Since $f_{z_k}-c\in C_0(N)$, the set $A_k\coloneqq\{m\in N : |f_{z_k}(m)-c|\ge\varepsilon/2\}$ is finite, and hence $B_k=\varphi_k^{-1}(A_k)$ is finite for every $k\ge 1$. Because $N$ is uncountable, we may pick $n_0\in N\setminus\bigcup_{k=1}^{\infty} B_k$; then $|f(z_k n_0)-c|<\varepsilon/2$ for all $k\ge 1$, whereas \eqref{eq:translate} applied to $n_0$ gives $|f(z_k n_0)-f(z_k)|<\varepsilon/2$ for all sufficiently large $k$, contradicting $|f(z_k)-c|\ge\varepsilon$. Thus $f-c\in C_0(G)$.
\end{proof}

\begin{example}[Sharpness of ascendancy]\label{ex:free_product}
Ascendancy cannot be dropped in Theorem~\ref{thm:ascendant}. Let $H$ be any nontrivial group and set $G\coloneqq\langle a\rangle * H$, where $\langle a\rangle\cong\mathbb{Z}$. Every element $x\in G\setminus\{e\}$ has a unique reduced normal form $x=x_1\cdots x_s$ ($s\ge 1$) with syllables alternating between $\langle a\rangle\setminus\{e\}$ and $H\setminus\{e\}$. Let $A\subseteq G$ be the set of reduced words whose first syllable is $a^k$ with $k\ge 1$. Fix $g\in G\setminus\{e\}$ with reduced form $g=y_1\cdots y_r$. If $x=x_1\cdots x_s\in A$ satisfies $xg\notin A$, then in forming $xg$ the entire suffix $x_2\cdots x_s$ (which is empty when $s=1$) must cancel against the prefix $y_1\cdots y_{s-1}$ and the syllable $y_s$ must cancel or reverse $x_1=a^k$; thus $1\le s\le r$, $x_2\cdots x_s=(y_1\cdots y_{s-1})^{-1}$, $y_s=a^{-m}$ for some $m\ge 1$, and $1\le k\le m$. Hence $\{x\in A : xg\notin A\}$ is finite for every $g\in G$. Applying this to $g$ and $g^{-1}$ shows that $Ag\triangle A$ is finite for every $g\in G$, so $\chi_A\in C_h(G)$. Since both $A$ (containing $\{a^k : k\ge 1\}$) and $G\setminus A$ (containing $\{a^{-k} : k\ge 1\}$) are infinite, $[\chi_A]$ is a nontrivial projection in $C(\partial G)$, so $\partial G$ is disconnected.

Taking $H$ to be an uncountable abelian group gives a group $G=\langle a\rangle * H$ containing an uncountable subgroup $H$ with singleton Higson corona while $\partial G$ is disconnected. Taking $H=F_\kappa$ (for any cardinal $\kappa\ge 1$) shows that every free group of rank at least $2$ has a disconnected Higson corona.
\end{example}

\begin{example}[Non-permanence under quotients by countably infinite normal subgroups]\label{ex:quotient}
While quotienting an infinite group by a finite normal subgroup is a coarse equivalence and therefore preserves $\partial G$, quotients by countably infinite normal subgroups need not preserve singleton coronas. Let $G\coloneqq\mathbb{Z}\times F_{\aleph_1}$ and $N\coloneqq\mathbb{Z}\times\{e\}\trianglelefteq G$. The known infinite-center result, also recovered in Theorem~\ref{thm:commuting_product}, makes $\partial G$ a singleton, whereas $G/N\cong F_{\aleph_1}$ has a disconnected Higson corona by Example~\ref{ex:free_product}.
\end{example}

Recall that a group $G$ is \emph{hyperabelian} if it admits an ascending series of normal subgroups of $G$,
\[
  \{e\}=G_0\trianglelefteq G_1\trianglelefteq\cdots\trianglelefteq G_\alpha\trianglelefteq\cdots\trianglelefteq G_\lambda=G,
\]
such that $G_\beta=\bigcup_{\alpha<\beta} G_\alpha$ for every limit ordinal $\beta\le\lambda$ and $G_{\alpha+1}/G_\alpha$ is abelian for every $\alpha<\lambda$.

\begin{corollary}\label{cor:hyperabelian}
Let $G$ be an uncountable group.
\begin{enumerate}
  \item If $G$ is virtually hyperabelian (in particular, virtually solvable), then $\partial G$ is a singleton.
  \item If $G$ is periodic and contains an involution with finite centralizer, then $\partial G$ is a singleton.
\end{enumerate}
\end{corollary}

\begin{proof}
For (1), by Proposition~\ref{prop:finite_index} we may pass to an uncountable finite-index subgroup and assume that $G$ is hyperabelian, with ascending normal series $(G_\alpha)_{\alpha\le\lambda}$. Let $\beta\le\lambda$ be the least ordinal such that $G_\beta$ is uncountable.
\begin{itemize}
  \item If $\beta=\gamma+1$ is a successor ordinal, then $G_\gamma$ is countable and $G_\beta/G_\gamma$ is abelian. For every countable subset $S\subseteq G_\beta$, the subgroup $H\coloneqq\langle S,G_\gamma\rangle\le G_\beta$ is countable and normal in $G_\beta$, so $G_\beta$ is $\omega_1$-normal.
  \item If $\beta$ is a limit ordinal, then $G_\beta=\bigcup_{\alpha<\beta} G_\alpha$ is uncountable while each $G_\alpha$ ($\alpha<\beta$) is countable, so $\beta$ has uncountable cofinality. Thus every countable subset $S\subseteq G_\beta$ lies in some countable normal subgroup $G_\alpha\trianglelefteq G_\beta$ ($\alpha<\beta$), so $G_\beta$ is again $\omega_1$-normal.
\end{itemize}
In either case, $\partial G_\beta$ is a singleton by Corollary~\ref{cor:transporter_apps}(2). Since $G_\beta\trianglelefteq G$ is uncountable, Theorem~\ref{thm:ascendant} shows that $\partial G$ is a singleton.

For (2), by \v{S}unkov's theorem \citep{Shunkov1972}, every periodic group containing an involution with finite centralizer is locally finite and virtually solvable; apply (1).
\end{proof}

\begin{example}[Separation of Corollary~\ref{cor:hyperabelian} from Theorem~\ref{thm:transporter}]\label{ex:metabelian}
Let $I$ be an uncountable set, put $V\coloneqq\bigoplus_{i\in I}\mathbb{F}_4$, and let $C_3=\langle\zeta\rangle$ act on $V$ by scalar multiplication by a primitive cube root of unity $\zeta\in\mathbb{F}_4^\times$. Set $G\coloneqq V\rtimes C_3$, with multiplication $(v,\zeta^r)(w,\zeta^s)=(v+\zeta^r w,\zeta^{r+s})$. Then $G$ is uncountable, locally finite, and metabelian, so $\partial G$ is a singleton by Corollary~\ref{cor:hyperabelian}(1). However, $G$ fails the hypothesis of Theorem~\ref{thm:transporter}: for $S\coloneqq\{(0,\zeta)\}$ and any $t=(v,\zeta^r)\in G$, we have $t^{-1}=(-\zeta^{-r}v,\zeta^{-r})$ and hence
\[
  t^{-1}(0,\zeta)t = \bigl(\zeta^{-r}(\zeta-1)v,\zeta\bigr).
\]
Since $\zeta-1\neq 0$ in $\mathbb{F}_4$, the map $v\mapsto\zeta^{-r}(\zeta-1)v$ is a bijection of $V$ for each $r\in\{0,1,2\}$. Hence $|T_G(S,H)|\le 3|H|\le\aleph_0$ for every countable subset $H\subseteq G$.
\end{example}

\section{Coset connectivity and groups of cardinality at least \texorpdfstring{$\aleph_2$}{aleph-2}}\label{sec:aleph2}

We use the bipartite coset-graph formulation of Holt's periodic argument \citep[\S3]{Holt1981}, as developed by \citet[Lemma~10.1]{Cornulier} and explained by \citet{RashedHolt}. The graph construction and alternating path below are shared ingredients; Theorem~\ref{thm:nilpotent_cosets} adds the nilpotent case by induction through a central quotient.

Let $\Gamma$ be a group generated by two proper subgroups $H$ and $K$. Consider the \emph{bipartite coset graph} $\mathcal{X}(\Gamma;H,K)$ whose vertex set is the disjoint union
\[
  V\bigl(\mathcal{X}(\Gamma;H,K)\bigr) \coloneqq H\backslash\Gamma \;\sqcup\; K\backslash\Gamma
\]
of right cosets, and whose edge set is $(H\cap K)\backslash\Gamma$, with the edge $(H\cap K)g$ joining $Hg$ to $Kg$. Since $(H\cap K)g_1=(H\cap K)g_2$ holds if and only if both $Hg_1=Hg_2$ and $Kg_1=Kg_2$, this defines a simple bipartite graph with no isolated vertices. Right multiplication by $\Gamma$ acts on $\mathcal{X}(\Gamma;H,K)$ by graph automorphisms and is transitive on edges; the connected component containing the edge $H\cap K$ is preserved by both $H$ and $K$, hence by $\Gamma=\langle H,K\rangle$, so $\mathcal{X}(\Gamma;H,K)$ is connected.

Let $\mathcal{X}^\circ(\Gamma;H,K)$ denote the induced subgraph obtained by deleting the two base vertices $H$ and $K$ and all edges incident to them. Since $\mathcal{X}(\Gamma;H,K)$ is connected, every connected component of $\mathcal{X}^\circ(\Gamma;H,K)$ contains at least one vertex adjacent in $\mathcal{X}(\Gamma;H,K)$ to $H$ or $K$; these neighbor sets are
\[
  N_H^\circ \coloneqq \{Kh : h\in H\setminus K\}
  \qquad\text{and}\qquad
  N_K^\circ \coloneqq \{Hk : k\in K\setminus H\},
\]
both nonempty because $H$ and $K$ are proper subgroups generating $\Gamma$. Consequently, $\mathcal{X}^\circ(\Gamma;H,K)$ is connected if and only if every $Hk\in N_K^\circ$ is joined to every $Kh\in N_H^\circ$ by a walk (equivalently, a path) in $\mathcal{X}^\circ(\Gamma;H,K)$.

We need two lemmas. The first is a special case of B.~H.~Neumann's finite coset-covering theorem \citep{Neumann1954}; its short proof is included.

\begin{lemma}[Two-subgroup covering lemma]\label{lem:neumann}
Let $\Delta$ be a group and let $A,B\le\Delta$ be subgroups of infinite index. Then $\Delta$ cannot be covered by finitely many left cosets of $A$ together with finitely many left cosets of $B$.
\end{lemma}

\begin{proof}
Suppose $\Delta=\bigcup_{i=1}^r a_i A\cup\bigcup_{j=1}^s b_j B$ with $r,s\ge 0$. Since $[\Delta:A]=\infty$, we may choose $x\in\Delta\setminus\bigcup_{i=1}^r a_i A$. Because left cosets of $A$ are pairwise disjoint, $xA\cap\bigcup_{i=1}^r a_i A=\varnothing$, so $xA\subseteq\bigcup_{j=1}^s b_j B$. Multiplying on the left by $a_i x^{-1}$ gives $a_i A\subseteq\bigcup_{j=1}^s (a_i x^{-1}b_j)B$ for each $1\le i\le r$. Hence $\Delta$ is covered by $(r+1)s$ left cosets of $B$, contradicting $[\Delta:B]=\infty$.
\end{proof}

The second lemma isolates the alternating-path step in Holt's proof. The same path requires only one of the two cyclic coset orbits to be finite, which is the form needed below.

\begin{lemma}[One-sided Holt zigzag]\label{lem:one_sided_holt}
Let $\Delta$ be a group generated by two proper subgroups $H$ and $K$. Fix $h\in H\setminus K$ and $k\in K\setminus H$, and set $s\coloneqq h^{-1}k$. If at least one of the right $\langle s\rangle$-orbits of $H$ or $K$ is finite, then $Hk$ and $Kh$ are joined by a path in $\mathcal{X}^\circ(\Delta;H,K)$.
\end{lemma}

\begin{proof}
Notice that $s\notin H$ (since $s\in H$ gives $k=hs\in H$) and $s\notin K$ (since $s\in K$ gives $h^{-1}=sk^{-1}\in K$). Interchanging $(H,h)$ and $(K,k)$ replaces $s$ by $k^{-1}h=s^{-1}$, leaves $\langle s\rangle$ unchanged, and swaps $Hk$ with $Kh$; thus we may assume that the right $\langle s\rangle$-orbit of $H$ is finite. Let $p\ge 2$ be the least positive integer such that $Hs^p=H$, and let $q\in\{2,3,\ldots\}\cup\{\infty\}$ be the least positive integer with $Ks^q=K$ (putting $q\coloneqq\infty$ if $Ks^i\neq K$ for all $i\ge 1$).

If $q\le p$, consider the vertex sequence
\[
  Hs,\ Ks,\ Hs^2,\ Ks^2,\ \ldots,\ Hs^{q-1},\ Ks^{q-1}.
\]
For $1\le i\le q-1$, minimality of $p$ and $q$ ensures $Hs^i\neq H$ and $Ks^i\neq K$, so all vertices lie in $\mathcal{X}^\circ(\Delta;H,K)$. Consecutive vertices are joined by the edges $(H\cap K)s^i$ ($1\le i\le q-1$) and $(H\cap K)ks^i$ ($1\le i\le q-2$, joining $Kks^i=Ks^i$ to $Hks^i=Hhs^{i+1}=Hs^{i+1}$), and the endpoints are $Hs=Hh^{-1}k=Hk$ and $Ks^{q-1}=Ks^{-1}=Kk^{-1}h=Kh$.

If $p<q$ (including $q=\infty$), consider instead the vertex sequence
\[
  Ks^{-1},\ Hs^{-1},\ Ks^{-2},\ Hs^{-2},\ \ldots,\ Ks^{-(p-1)},\ Hs^{-(p-1)}.
\]
For $1\le i\le p-1$, we have $Hs^{-i}\neq H$ (since $Hs^{-i}=H\iff Hs^i=H$) and $Ks^{-i}\neq K$ (since $Ks^{-i}=K\iff Ks^i=K$), so all vertices lie in $\mathcal{X}^\circ(\Delta;H,K)$. Consecutive vertices are joined by the edges $(H\cap K)s^{-i}$ ($1\le i\le p-1$) and $(H\cap K)hs^{-i}$ ($1\le i\le p-2$, joining $Hhs^{-i}=Hs^{-i}$ to $Khs^{-i}=K(ks^{-1})s^{-i}=Ks^{-(i+1)}$). Finally, the endpoints are $Ks^{-1}=Kh$ and $Hs^{-(p-1)}=Hs^p s^{-(p-1)}=Hs=Hk$.
\end{proof}

\begin{theorem}[Nilpotent coset-connectivity theorem]\label{thm:nilpotent_cosets}
Let $\Delta$ be a nilpotent group generated by two proper subgroups $H$ and $K$. Then the deleted bipartite coset graph $\mathcal{X}^\circ(\Delta;H,K)$ is connected.
\end{theorem}

\begin{proof}
We proceed by induction on the nilpotency class $c\ge 1$ of $\Delta$ (note $c\ge 1$ since $H,K<\Delta$ are proper).

If $c=1$, then $\Delta$ is abelian. For any $h\in H\setminus K$ and $k\in K\setminus H$, the edge $(H\cap K)hk$ joins $H(hk)=Hk\neq H$ to $K(hk)=K(kh)=Kh\neq K$, so $\mathcal{X}^\circ(\Delta;H,K)$ is connected.

Now assume $c\ge 2$ and that the result holds for all nilpotent groups of class $c-1$. Define the lower central series inductively by $\gamma_1(\Delta)\coloneqq\Delta$ and $\gamma_{i+1}(\Delta)\coloneqq[\gamma_i(\Delta),\Delta]$, and put $Z_0\coloneqq\gamma_c(\Delta)$. Then $\{e\}\neq Z_0\le Z(\Delta)$, and $\gamma_{c-1}(\Delta)\neq Z_0$ (otherwise $Z_0=[\gamma_{c-1}(\Delta),\Delta]=[Z_0,\Delta]=\gamma_{c+1}(\Delta)=\{e\}$, contradicting class $c$). Thus $\bar{\Delta}\coloneqq\Delta/Z_0$ has nilpotency class $c-1$. Set $\bar{H}\coloneqq HZ_0/Z_0$ and $\bar{K}\coloneqq KZ_0/Z_0$, writing $\bar{x}\coloneqq xZ_0$ for $x\in\Delta$. We consider two cases.

\emph{Case 1: $[\Delta:HZ_0]<\infty$ (or symmetrically $[\Delta:KZ_0]<\infty$).}

\emph{Claim 1. Every element of $\Delta$ has a positive power in $H$.}
Let $L\coloneqq\bigcap_{g\in\Delta} g(HZ_0)g^{-1}\trianglelefteq\Delta$ be the kernel of the left translation action of $\Delta$ on $\Delta/HZ_0$; then $m\coloneqq[\Delta:L]\le[\Delta:HZ_0]!<\infty$, so $g^m\in L\subseteq HZ_0$ for every $g\in\Delta$.

For each $i\ge 1$, the inclusion $[\gamma_{i+1}(\Delta),\gamma_{i+1}(\Delta)]\le[\gamma_{i+1}(\Delta),\Delta]=\gamma_{i+2}(\Delta)$ shows that $\gamma_{i+1}(\Delta)/\gamma_{i+2}(\Delta)$ is abelian. Recall the commutator identities $[xy,z]=x[y,z]x^{-1}[x,z]$ and $[x,yz]=[x,y]\,y[x,z]y^{-1}$. For $x,y\in\gamma_i(\Delta)$ and $z\in\Delta$, we have $[y,z]\in\gamma_{i+1}(\Delta)$, so $x[y,z]x^{-1}=[x,[y,z]][y,z]\equiv[y,z]\pmod{\gamma_{i+2}(\Delta)}$ and hence $[xy,z]\equiv[y,z][x,z]\equiv[x,z][y,z]\pmod{\gamma_{i+2}(\Delta)}$; similarly, $y[x,z]y^{-1}\equiv[x,z]\pmod{\gamma_{i+2}(\Delta)}$ for $x\in\gamma_i(\Delta)$ and $y,z\in\Delta$, giving $[x,yz]\equiv[x,y][x,z]\pmod{\gamma_{i+2}(\Delta)}$. Moreover, $[w,z]\in\gamma_{i+2}(\Delta)$ for every $w\in\gamma_{i+1}(\Delta)$ and $z\in\Delta$. Hence $(x\gamma_{i+1}(\Delta),z)\mapsto[x,z]\gamma_{i+2}(\Delta)$ is a well-defined bihomomorphism $(\gamma_i(\Delta)/\gamma_{i+1}(\Delta))\times\Delta\to\gamma_{i+1}(\Delta)/\gamma_{i+2}(\Delta)$ whose image generates $\gamma_{i+1}(\Delta)/\gamma_{i+2}(\Delta)$. Writing left-normed simple commutators inductively as $[g_1]\coloneqq g_1$ and $[g_1,\ldots,g_i]\coloneqq[[g_1,\ldots,g_{i-1}],g_i]$, induction on $i$ from $1$ to $c$ (where $\gamma_{c+1}(\Delta)=\{e\}$) shows that $(g_1,\ldots,g_c)\mapsto[g_1,\ldots,g_c]$ is a homomorphism in each coordinate from $\Delta^c$ into the abelian group $Z_0$ whose image generates $Z_0$.

For any $g_1,\ldots,g_c\in\Delta$, write $g_j^m=h_j z_j$ with $h_j\in H$ and $z_j\in Z_0\le Z(\Delta)$. Since $c\ge 2$ and multiplying any entry of a commutator by a central element leaves the commutator unchanged,
\[
  [g_1,\ldots,g_c]^{m^c}
  = [g_1^m,\ldots,g_c^m]
  = [h_1 z_1,\ldots,h_c z_c]
  = [h_1,\ldots,h_c]\in H\cap Z_0.
\]
Because $Z_0$ is abelian and generated by the simple $c$-fold commutators $[g_1,\ldots,g_c]$, it follows that $z^{m^c}\in H\cap Z_0$ for every $z\in Z_0$. For any $s\in\Delta$, write $s^m=h_0z$ with $h_0\in H$ and $z\in Z_0$. Since $z$ is central,
\[
  s^{m^{c+1}} = (h_0 z)^{m^c} = h_0^{m^c} z^{m^c}\in H,
\]
This proves Claim~1. In particular, for $h\in H\setminus K$, $k\in K\setminus H$, and $s=h^{-1}k$, the right $\langle s\rangle$-orbit of $H$ is finite. Lemma~\ref{lem:one_sided_holt} joins $Hk$ to $Kh$ in $\mathcal{X}^\circ(\Delta;H,K)$, proving connectivity in Case~1.

\emph{Case 2: $[\Delta:HZ_0]=\infty$ and $[\Delta:KZ_0]=\infty$.}
Then $\bar{H}$ and $\bar{K}$ are proper subgroups generating $\bar{\Delta}$, so $\mathcal{X}^\circ(\bar{\Delta};\bar{H},\bar{K})$ is connected by the induction hypothesis. Let $\mathcal{X}^{\circ\circ}$ be the induced subgraph of $\mathcal{X}^\circ(\Delta;H,K)$ on
\[
  V(\mathcal{X}^{\circ\circ})
  \coloneqq \{Hx : x\in\Delta\setminus HZ_0\} \;\sqcup\; \{Kx : x\in\Delta\setminus KZ_0\}.
\]
The projection $\pi\colon\mathcal{X}^{\circ\circ}\to\mathcal{X}^\circ(\bar{\Delta};\bar{H},\bar{K})$ defined by $\pi(Hx)\coloneqq\bar{H}\bar{x}$ and $\pi(Kx)\coloneqq\bar{K}\bar{x}$ is a graph homomorphism preserving the bipartition. For each $z\in Z_0$, right multiplication $v\mapsto v\cdot z$ restricts to a graph automorphism of $\mathcal{X}^{\circ\circ}$ satisfying $\pi(v\cdot z)=\pi(v)$.

\emph{Claim 2. Every path in the quotient graph lifts from each preimage of its initial vertex.}
Let $v=Hx\in V(\mathcal{X}^{\circ\circ})$, and let $(\bar{H}\cap\bar{K})\bar{g}$ be an edge of $\mathcal{X}^\circ(\bar{\Delta};\bar{H},\bar{K})$ incident to $\bar{H}\bar{x}=\bar{H}\bar{g}$. Then $g\in HZ_0 x=Z_0 Hx$, so we may write $g=z_1 h_1 x$ with $z_1\in Z_0$ and $h_1\in H$. Put $g'\coloneqq z_1^{-1}g=h_1 x$. Then $Hg'=Hx=v$ and $\bar{g}'=\bar{g}$, so $\bar{K}\bar{g}'=\bar{K}\bar{g}\neq\bar{K}$, whence $g'\notin KZ_0$. Thus the edge $(H\cap K)g'$ of $\mathcal{X}^{\circ\circ}$ joins $v$ to $Kg'\in V(\mathcal{X}^{\circ\circ})$ over $(\bar{H}\cap\bar{K})\bar{g}$. Symmetrically, edges lift from vertices of the form $Kx\in V(\mathcal{X}^{\circ\circ})$. By induction on path length, every path in $\mathcal{X}^\circ(\bar{\Delta};\bar{H},\bar{K})$ lifts to $\mathcal{X}^{\circ\circ}$ starting at any preimage of its initial vertex.

\emph{Claim 3. The graph $\mathcal{X}^{\circ\circ}$ is connected.}
We first connect the vertices in one central fiber and then use Claim~2 to reach that fiber from every other vertex. By Lemma~\ref{lem:neumann}, choose $g_0\in\Delta\setminus(HZ_0\cup KZ_0)$; then $Kg_0\in V(\mathcal{X}^{\circ\circ})$. Fix any $z\in Z_0$. Because $\mathcal{X}(\Delta;H,K)$ is connected, choose a path $\gamma$ in $\mathcal{X}(\Delta;H,K)$ from $K$ to $Kz$ with consecutive edge representatives $u_1,\ldots,u_t\in\Delta$ ($t\ge 0$). By Lemma~\ref{lem:neumann} applied to $HZ_0$ and $KZ_0$, we may choose
\[
  g\in\Delta\setminus\left(HZ_0\cup KZ_0\cup\bigcup_{i=1}^t\bigl(u_i^{-1}HZ_0\cup u_i^{-1}KZ_0\bigr)\right).
\]
Then $g,zg\notin KZ_0$ and $u_i g\notin HZ_0\cup KZ_0$ for all $1\le i\le t$, so the right-translated path $\gamma\cdot g$ lies entirely in $\mathcal{X}^{\circ\circ}$ and connects $Kg$ to $Kzg$. Since $\mathcal{X}^\circ(\bar{\Delta};\bar{H},\bar{K})$ is connected, pick a path from $\bar{K}\bar{g}$ to $\bar{K}\bar{g}_0$ and lift it from $Kg$ to a path $\sigma$ in $\mathcal{X}^{\circ\circ}$. Because $\pi$ preserves the bipartition, the endpoint of $\sigma$ is a $K$-coset $Ky$ with $\bar{K}\bar{y}=\bar{K}\bar{g}_0$, so $y\in KZ_0 g_0$ and hence $Ky=Kz_1 g_0$ for some $z_1\in Z_0$. Right-multiplying $\sigma$ by $z\in Z_0\le Z(\Delta)$ produces a path $\sigma\cdot z$ in $\mathcal{X}^{\circ\circ}$ from $Kzg$ to $Kz_1 g_0 z=Kzz_1 g_0$. Concatenating $\sigma^{-1}$, $\gamma\cdot g$, and $\sigma\cdot z$, and then right-multiplying the resulting walk by $z_1^{-1}\in Z_0\le Z(\Delta)$ (so that $Kz_1 g_0 z_1^{-1}=Kg_0$ and $Kzz_1 g_0 z_1^{-1}=Kzg_0$), yields a walk---and hence a path---in $\mathcal{X}^{\circ\circ}$ from $Kg_0$ to $Kzg_0$. Finally, for any $v\in V(\mathcal{X}^{\circ\circ})$, lifting a path in $\mathcal{X}^\circ(\bar{\Delta};\bar{H},\bar{K})$ from $\pi(v)$ to $\bar{K}\bar{g}_0$ connects $v$ to $Kzg_0$ for some $z\in Z_0$, and hence to $Kg_0$. Thus $\mathcal{X}^{\circ\circ}$ is connected.

\emph{Attaching the remaining vertices.}
It remains to connect every vertex in $V(\mathcal{X}^\circ(\Delta;H,K))\setminus V(\mathcal{X}^{\circ\circ})$ to $\mathcal{X}^{\circ\circ}$. Every such vertex has the form $Hz$ ($z\in Z_0\setminus H$) or $Kz$ ($z\in Z_0\setminus K$). Since $\bar{H}$ and $\bar{K}$ are proper subgroups generating $\bar{\Delta}$, neither contains the other; in particular, $\bar{H}\not\subseteq\bar{K}$ implies $H\not\subseteq KZ_0$ (otherwise $HZ_0\subseteq KZ_0$), and symmetrically $K\not\subseteq HZ_0$. Choose $h_1\in H\setminus KZ_0$ and $k_1\in K\setminus HZ_0$. For $z\in Z_0\setminus H$, we have $h_1 z\notin KZ_0$, so the edge $(H\cap K)(h_1 z)$ joins $Hz$ to $K(h_1 z)\in V(\mathcal{X}^{\circ\circ})$. Symmetrically, for $z\in Z_0\setminus K$, the edge $(H\cap K)(k_1 z)$ joins $Kz$ to $H(k_1 z)\in V(\mathcal{X}^{\circ\circ})$. Therefore $\mathcal{X}^\circ(\Delta;H,K)$ is connected.
\end{proof}

\begin{corollary}[Generalized coset-connectivity lemma]\label{cor:generalized_holt}
Let $\Gamma$ be a group in which every $2$-generated subgroup is periodic or nilpotent (in particular, $\Gamma$ may be any periodic group or any locally nilpotent group), and suppose $\Gamma=\langle H,K\rangle$ for two proper subgroups $H$ and $K$. Then $\mathcal{X}^\circ(\Gamma;H,K)$ is connected. Consequently, if $S$ is a set and $\alpha\colon\Gamma\to S$ is constant on each right coset of $H$ other than $H$ itself, and on each right coset of $K$ other than $K$ itself, then $\alpha$ is constant on $\Gamma\setminus(H\cap K)$; if $\alpha$ is also constant on at least one of $H$ and $K$, then it is constant on $\Gamma$.
\end{corollary}

\begin{proof}
Fix $h\in H\setminus K$ and $k\in K\setminus H$, and set $\Delta\coloneqq\langle h,k\rangle$. If $\Delta$ is periodic, then $s\coloneqq h^{-1}k\in\Delta$ has finite order, so Lemma~\ref{lem:one_sided_holt} joins $Hk$ to $Kh$ in $\mathcal{X}^\circ(\Gamma;H,K)$. If $\Delta$ is nilpotent, put $H_0\coloneqq H\cap\Delta$ and $K_0\coloneqq K\cap\Delta$; then $H_0$ and $K_0$ are proper subgroups generating $\Delta$ (since $h\in H_0\setminus K_0$ and $k\in K_0\setminus H_0$). For $g_1,g_2\in\Delta$, we have $H_0 g_1=H_0 g_2\iff g_1 g_2^{-1}\in H\cap\Delta\iff Hg_1=Hg_2$, and similarly for $K_0$ and $H_0\cap K_0=H\cap K\cap\Delta$. Hence the natural map $\iota\colon\mathcal{X}(\Delta;H_0,K_0)\to\mathcal{X}(\Gamma;H,K)$ sending $H_0 g\mapsto Hg$, $K_0 g\mapsto Kg$, and $(H_0\cap K_0)g\mapsto(H\cap K)g$ is a well-defined graph homomorphism that is injective on vertices and maps $\{H_0,K_0\}$ to $\{H,K\}$. In particular, $\iota$ maps $\mathcal{X}^\circ(\Delta;H_0,K_0)$ into $\mathcal{X}^\circ(\Gamma;H,K)$. By Theorem~\ref{thm:nilpotent_cosets}, $H_0 k$ and $K_0 h$ are joined in $\mathcal{X}^\circ(\Delta;H_0,K_0)$, so $Hk$ and $Kh$ are joined in $\mathcal{X}^\circ(\Gamma;H,K)$. Thus $\mathcal{X}^\circ(\Gamma;H,K)$ is connected.

For the second assertion, $\alpha$ induces a well-defined map $\bar{\alpha}\colon V(\mathcal{X}^\circ(\Gamma;H,K))\to S$ with $\bar{\alpha}(Hg)=\alpha(g)$ ($Hg\neq H$) and $\bar{\alpha}(Kg)=\alpha(g)$ ($Kg\neq K$). For every edge $(H\cap K)g$ of $\mathcal{X}^\circ(\Gamma;H,K)$, we have $\bar{\alpha}(Hg)=\alpha(g)=\bar{\alpha}(Kg)$. Since $\mathcal{X}^\circ(\Gamma;H,K)$ is connected and nonempty (containing $N_H^\circ$ and $N_K^\circ$), $\bar{\alpha}$ takes a single constant value $c\in S$. Every $g\in\Gamma\setminus(H\cap K)$ satisfies $Hg\neq H$ or $Kg\neq K$, giving $\alpha(g)=c$. Finally, if $\alpha$ is constant on $K$ (or symmetrically on $H$), picking $k_0\in K\setminus H\subseteq\Gamma\setminus(H\cap K)$ shows that the constant value on $K$ is $\alpha(k_0)=c$, whence $\alpha\equiv c$ on $\Gamma$.
\end{proof}

\begin{remark}\label{rem:dihedral_sharp}
Nilpotency cannot be weakened to virtually abelian or metabelian in Theorem~\ref{thm:nilpotent_cosets}: in the infinite dihedral group $D_\infty=\langle h,k \mid h^2=k^2=e\rangle$, the bipartite coset graph $\mathcal{X}(D_\infty;\langle h\rangle,\langle k\rangle)$ is a bi-infinite line, so deleting the adjacent vertices $\langle h\rangle$ and $\langle k\rangle$ disconnects it into two rays.
\end{remark}

Because Corollary~\ref{cor:generalized_holt} does not require $H$ or $K$ to be finite, applying it to countable defect hulls forces any two incomparable defect hulls to collapse a left-slowly oscillating function.

\begin{lemma}[Incomparable defect hulls]\label{lem:incomparable}
Let $G$ be a group in which every $2$-generated subgroup is periodic or nilpotent, let $\alpha\colon G\to\mathbb{C}$ be left-slowly oscillating, and let $H,K\le G$ be defect hulls for $\alpha$ that are incomparable under inclusion ($H\not\subseteq K$ and $K\not\subseteq H$). Put $\Gamma\coloneqq\langle H,K\rangle$. Then $\alpha|_\Gamma$ is constant modulo $C_0(\Gamma)$.
\end{lemma}

\begin{proof}
Since $H$ and $K$ are incomparable, both are proper subgroups of $\Gamma$. Every $2$-generated subgroup of $\Gamma$ is periodic or nilpotent, and because $H$ and $K$ are defect hulls, $\alpha|_\Gamma$ is constant on every right coset of $H$ other than $H$ and on every right coset of $K$ other than $K$. Corollary~\ref{cor:generalized_holt} yields a constant $c\in\mathbb{C}$ such that $\alpha(x)=c$ for all $x\in\Gamma\setminus(H\cap K)$.

Choose $g\in H\setminus K$. For every $x\in H\cap K$, we have $gx\in H\setminus K\subseteq\Gamma\setminus(H\cap K)$, so $\alpha(gx)=c$. Consequently, for every $\varepsilon>0$,
\[
  \{x\in\Gamma : |\alpha(x)-c|\ge\varepsilon\}
  = \{x\in H\cap K : |\alpha(x)-\alpha(gx)|\ge\varepsilon\}
  \subseteq D_{g,\varepsilon}(\alpha),
\]
which is finite. Hence $\alpha|_\Gamma-c\in C_0(\Gamma)$.
\end{proof}

\begin{theorem}\label{thm:aleph2}
Let $G$ be a group of cardinality $|G|\ge\aleph_2$ in which every $2$-generated subgroup is periodic or nilpotent (in particular, $G$ may be any periodic group or any locally nilpotent group of cardinality at least $\aleph_2$), equipped with the finitary coarse structure. Then the Higson corona $\partial G$ is a singleton.
\end{theorem}

\begin{proof}
Fix $f\in C_h(G)$ and set $\alpha(x)\coloneqq f(x^{-1})$. By Lemma~\ref{lem:inversion}, $\alpha$ is left-slowly oscillating and it suffices to show that $\alpha$ is constant modulo $C_0(G)$. Suppose for contradiction that it is not. By Lemma~\ref{lem:oscillation}(1), there exist $\varepsilon>0$ and sequences $(u_n)_{n\ge 1}$ and $(v_n)_{n\ge 1}$ in $G$ with pairwise distinct combined terms such that $|\alpha(u_n)-\alpha(v_n)|\ge\varepsilon$ for all $n\ge 1$. Set $S\coloneqq\{u_n,v_n : n\ge 1\}$ and $H\coloneqq\operatorname{cl}_\alpha(S)$. By Lemma~\ref{lem:defect_hull}, $H$ is a countable defect hull containing $S$, and $\alpha|_H$ is not constant modulo $C_0(H)$ (if $\alpha|_H-c\in C_0(H)$, then $\alpha(u_n)\to c$ and $\alpha(v_n)\to c$, contradicting $|\alpha(u_n)-\alpha(v_n)|\ge\varepsilon$).

For each $x\in G\setminus H$, put $K(x)\coloneqq\operatorname{cl}_\alpha(\{x\})$, which is a countable defect hull by Lemma~\ref{lem:defect_hull}. Since $x\in K(x)\setminus H$, we have $K(x)\not\subseteq H$. If $H\not\subseteq K(x)$, then Lemma~\ref{lem:incomparable} would imply that $\alpha|_{\langle H,K(x)\rangle}$---and hence $\alpha|_H$---is constant modulo $C_0$, a contradiction. Therefore $H\subseteq K(x)$ for every $x\in G\setminus H$.

In particular, $S\subseteq K(x)$ for every $x\in G\setminus H$, so $\alpha|_{K(x)}$ is not constant modulo $C_0(K(x))$. Applying Lemma~\ref{lem:incomparable} to any two subgroups in $\mathcal{L}\coloneqq\{K(x) : x\in G\setminus H\}$ shows that no two members of $\mathcal{L}$ are incomparable; hence $\mathcal{L}$ is totally ordered by inclusion.

Construct elements $(x_\beta)_{\beta<\omega_1}$ in $G\setminus H$ by transfinite recursion: having chosen $x_\gamma$ for all $\gamma<\beta$ (where $\beta<\omega_1$), the union $H\cup\bigcup_{\gamma<\beta} K(x_\gamma)$ is countable, so since $|G|\ge\aleph_2$ we may choose $x_\beta\in G\setminus\bigl(H\cup\bigcup_{\gamma<\beta} K(x_\gamma)\bigr)$. For $\gamma<\beta<\omega_1$, we have $x_\beta\notin K(x_\gamma)$, so $K(x_\beta)\not\subseteq K(x_\gamma)$ and total ordering forces $K(x_\gamma)\subsetneq K(x_\beta)$. Thus the union $U\coloneqq\bigcup_{\beta<\omega_1} K(x_\beta)$ has cardinality $\aleph_1$.

Since $|G|\ge\aleph_2$, we may choose $y\in G\setminus U$. Because $H\subseteq K(x_0)\subseteq U$, we have $y\in G\setminus H$, so $K(y)\in\mathcal{L}$. For every $\beta<\omega_1$, $y\notin K(x_\beta)$ gives $K(y)\not\subseteq K(x_\beta)$, and hence $K(x_\beta)\subsetneq K(y)$ by total ordering. Therefore $U\subseteq K(y)$, contradicting $|U|=\aleph_1$ and $|K(y)|\le\aleph_0$.
\end{proof}

\begin{remark}\label{rem:eastridge_periodic}
For periodic groups, the conclusion is already implied by \citep[Theorem~7.9]{Eastridge2017}. Theorem~\ref{thm:nonperiodic_2gen} below strengthens the nonperiodic case to every uncountable cardinality. The role of this section is the geometric argument, especially Theorem~\ref{thm:nilpotent_cosets}, rather than a new periodic cardinality threshold.
\end{remark}

\begin{remark}\label{rem:periodicity_sharp}
The $2$-generator hypothesis in Theorem~\ref{thm:aleph2} cannot be omitted at any cardinality: by Example~\ref{ex:free_product}, for every uncountable cardinal $\kappa$, the free group $F_\kappa$ has a disconnected Higson corona.
\end{remark}

\section{Nonperiodic locally nilpotent groups and rigidity at \texorpdfstring{$\aleph_1$}{aleph-1}}\label{sec:aleph1}

Recall that a subgroup $M\le G$ is \emph{almost malnormal} in $G$ if $M\cap k^{-1}Mk$ is finite for every $k\in G\setminus M$.

\begin{lemma}[Almost-malnormality rigidity]\label{lem:almost_malnormal}
Let $G$ be a group in which every $2$-generated subgroup is periodic or nilpotent, let $\alpha\colon G\to\mathbb{C}$ be left-slowly oscillating, and let $M<G$ be a defect hull for $\alpha$. Suppose there exists $k\in G\setminus M$ such that either:
\begin{enumerate}
  \item $M\cap k^{-1}Mk$ is infinite and $k^d\in M$ for some integer $d\ge 1$ (in particular, this holds if $k$ has finite order and $M\cap k^{-1}Mk$ is infinite); or
  \item $J\coloneqq M\cap k^{-1}Mk\cap kMk^{-1}$ is infinite.
\end{enumerate}
Then $\alpha|_M$ is constant modulo $C_0(M)$.
\end{lemma}

\begin{proof}
Suppose for contradiction that $\alpha|_M$ is not constant modulo $C_0(M)$. Put $K\coloneqq\langle k\rangle$ and $\Gamma\coloneqq\langle M,K\rangle$. Since $k\notin M$, $M$ is a proper subgroup of $\Gamma$. If $M\le K$, then $K$ is abelian and hence $k\in N_G(M)\setminus M$, contradicting Lemma~\ref{lem:self_normalizing}(2). Thus $K$ is also a proper subgroup of $\Gamma$.

In both (1) and (2), the intersection $M\cap k^{-1}Mk$ is infinite. We first show that if $y\in G$ satisfies $y\notin M$ and $ky\notin M$, then $\alpha(ky)=\alpha(y)$. Indeed, for every $h\in M\cap k^{-1}Mk$ we have $h\in M$ and $khk^{-1}\in M$. Because $M$ is a defect hull, $\alpha(hy)=\alpha(y)$ and $\alpha(k(hy))=\alpha((khk^{-1})(ky))=\alpha(ky)$, whence $|\alpha(k(hy))-\alpha(hy)|=|\alpha(ky)-\alpha(y)|$. Since $h\mapsto hy$ is injective on the infinite set $M\cap k^{-1}Mk$ and $D_{k,\varepsilon}(\alpha)$ is finite for every $\varepsilon>0$, this forces $\alpha(ky)=\alpha(y)$.

Next, we claim that $\alpha$ is constant on $Kx\setminus M$ for every $x\in\Gamma$. Notice that $Kx\setminus M$ is nonempty (otherwise $x,kx\in M$ and hence $k\in M$). We consider three cases.
\begin{itemize}
  \item If $Kx\cap M=\varnothing$, then the preceding paragraph gives $\alpha(k^{j+1}x)=\alpha(k^j x)$ for all $j\in\mathbb{Z}$, so $\alpha$ is constant on $Kx$.
  \item If $Kx\cap M\neq\varnothing$ and $k^d\in M$ for some integer $d\ge 1$, choose $m\in Kx\cap M$ so that $Kx=Km$, and let $d\ge 2$ be the least positive integer with $k^d\in M$. Because $\{j\in\mathbb{Z} : k^j\in M\}$ is a nontrivial subgroup of $\mathbb{Z}$ with least positive element $d$, we have $k^j m\in M$ if and only if $k^j\in M$, which holds if and only if $d\mid j$. Hence every element of $Km\setminus M$ has the form $k^{qd+r}m$ with $q\in\mathbb{Z}$ and $1\le r\le d-1$. For $1\le r\le d-2$ (vacuously when $d=2$), neither $k^r m$ nor $k^{r+1}m$ lies in $M$, so $\alpha(k^{r+1}m)=\alpha(k^r m)$, giving $\alpha(k^r m)=\alpha(km)$ for all $1\le r\le d-1$. Since $k^{qd}\in M$ and $k^r m\notin M$, the defect-hull property of $M$ yields $\alpha(k^{qd+r}m)=\alpha(k^{qd}(k^r m))=\alpha(k^r m)=\alpha(km)$, so $\alpha$ is constant on $Km\setminus M$.
  \item If $Kx\cap M\neq\varnothing$ and $k^d\notin M$ for every integer $d\ge 1$, then (1) fails and therefore (2) holds, so $J$ is infinite. Write $Kx=Km$ with $m\in M$. Because $k^{-d}=(k^d)^{-1}\notin M$ for every $d\ge 1$, we have $k^j m\in M$ if and only if $j=0$ (and $k$ has infinite order). Thus $Km\setminus M=\{k^j m : j\ge 1\}\cup\{k^j m : j\le -1\}$. By the zero-defect step, $\alpha(k^{j+1}m)=\alpha(k^j m)$ for all $j\ge 1$ and for all $j\le -2$, giving $\alpha(k^j m)=\alpha(km)$ for all $j\ge 1$ and $\alpha(k^j m)=\alpha(k^{-1}m)$ for all $j\le -1$. Fix $\varepsilon>0$. Since $J$ is infinite, $h\mapsto hm$ is injective on $J$, and $D_{k,\varepsilon/2}(\alpha)\cup D_{k^{-1},\varepsilon/2}(\alpha)$ is finite, we may choose $h\in J$ with $hm\notin D_{k,\varepsilon/2}(\alpha)\cup D_{k^{-1},\varepsilon/2}(\alpha)$. Because $khk^{-1}\in M$, $k^{-1}hk\in M$, $km\notin M$, and $k^{-1}m\notin M$, the defect-hull property of $M$ yields
  \[
    \begin{aligned}
      \alpha(k(hm))
        &= \alpha\bigl((khk^{-1})(km)\bigr) = \alpha(km),\\
      \alpha(k^{-1}(hm))
        &= \alpha\bigl((k^{-1}hk)(k^{-1}m)\bigr) = \alpha(k^{-1}m).
    \end{aligned}
  \]
  Consequently,
  \[
    |\alpha(km)-\alpha(k^{-1}m)|
    \le |\alpha(k(hm))-\alpha(hm)| + |\alpha(hm)-\alpha(k^{-1}(hm))|
    < \varepsilon.
  \]
  Since $\varepsilon>0$ was arbitrary, $\alpha(km)=\alpha(k^{-1}m)$, so $\alpha$ is constant on $Km\setminus M$.
\end{itemize}

Now define $\beta\colon\Gamma\to\mathbb{C}$ by
\[
  \beta(x) \coloneqq
  \begin{cases}
    \alpha(x)  & \text{if } x\in\Gamma\setminus M,\\
    \alpha(kx) & \text{if } x\in M.
  \end{cases}
\]
On every right $M$-coset other than $M$, $\beta$ agrees with $\alpha$ and is therefore constant. We claim that $\beta$ is constant on every right $K$-coset $Kx\subseteq\Gamma$ (including $K$ itself). If $Kx\cap M=\varnothing$, then $\beta|_{Kx}=\alpha|_{Kx}$ is constant. If $Kx\cap M\neq\varnothing$, write $Kx=Km$ with $m\in M$. For $x\in Km\setminus M$, we have $\beta(x)=\alpha(x)=\alpha(km)$. For $x\in Km\cap M$, we have $kx\in Km\setminus M$ (since $x\in M$ and $k\notin M$), so $\beta(x)=\alpha(kx)=\alpha(km)$. Thus $\beta$ is constant on $Km$.

Applying Corollary~\ref{cor:generalized_holt} to $(\Gamma,M,K,\beta)$, we conclude that $\beta$ takes a single constant value $c\in\mathbb{C}$ on $\Gamma$. In particular, $\alpha(km)=\beta(m)=c$ for every $m\in M$. Hence, for every $\varepsilon>0$,
\[
  \{m\in M : |\alpha(m)-c|\ge\varepsilon\}
  = \{m\in M : |\alpha(m)-\alpha(km)|\ge\varepsilon\}
  \subseteq D_{k,\varepsilon}(\alpha),
\]
which is finite. Therefore $\alpha|_M-c\in C_0(M)$, a contradiction.
\end{proof}

Lemma~\ref{lem:almost_malnormal}(2) resolves all nonperiodic locally nilpotent groups at every uncountable cardinality, including $\aleph_1$.

\begin{theorem}[Nonperiodic groups with periodic or nilpotent $2$-generated subgroups]\label{thm:nonperiodic_2gen}
Let $G$ be an uncountable group in which every $2$-generated subgroup is periodic or nilpotent. If $G$ is not periodic, then the Higson corona $\partial G$ is a singleton. In particular, every uncountable nonperiodic (e.g., torsion-free) locally nilpotent group has a singleton Higson corona.
\end{theorem}

\begin{proof}
Fix an element $a\in G$ of infinite order, and suppose for contradiction that some $f\in C_h(G)$ is not constant modulo $C_0(G)$. Set $\alpha(x)\coloneqq f(x^{-1})$. By Lemmas~\ref{lem:inversion} and~\ref{lem:oscillation}(1), there exists a countable subset $S\subseteq G$ such that $\alpha|_S$ is not constant modulo $C_0(S)$. Put $M\coloneqq\operatorname{cl}_\alpha(S\cup\{a\})$. By Lemma~\ref{lem:defect_hull}, $M$ is a countable defect hull containing $S\cup\{a\}$. Since $G$ is uncountable, $M<G$ is proper, and since $S\subseteq M$, $\alpha|_M$ is not constant modulo $C_0(M)$.

Choose $g\in G\setminus M$ and set $\Delta\coloneqq\langle a,g\rangle\le G$. Because $a\in\Delta$ has infinite order, $\Delta$ is not periodic and is therefore nilpotent, of some class $c\ge 1$. Put $M_0\coloneqq M\cap\Delta$; then $\langle a\rangle\le M_0$ is infinite, while $g\in\Delta\setminus M_0$, so $M_0<\Delta$ is proper.

Consider the upper central series of $\Delta$, defined inductively by $Z_0(\Delta)\coloneqq\{e\}$ and $Z_{i+1}(\Delta)/Z_i(\Delta)\coloneqq Z(\Delta/Z_i(\Delta))$, with $Z_c(\Delta)=\Delta$. Since $Z_0(\Delta)\le M_0$ and $Z_c(\Delta)=\Delta\not\subseteq M_0$, let $r\in\{0,\ldots,c-1\}$ be maximal with $Z_r(\Delta)\le M_0$, and choose $k\in Z_{r+1}(\Delta)\setminus M_0$. Because $k\in\Delta\setminus(M\cap\Delta)$, we have $k\in G\setminus M$. Moreover, since $Z_{r+1}(\Delta)$ is a subgroup, both $k$ and $k^{-1}$ lie in $Z_{r+1}(\Delta)$, and $[Z_{r+1}(\Delta),\Delta]\le Z_r(\Delta)\le M_0$. Hence, for every $m\in M_0$,
\[
  kmk^{-1} = [k,m]m \in Z_r(\Delta)M_0 = M_0
  \qquad\text{and}\qquad
  k^{-1}mk = [k^{-1},m]m \in Z_r(\Delta)M_0 = M_0.
\]
Thus $M_0=k^{-1}M_0k=kM_0k^{-1}$, which gives
\[
  \langle a\rangle \le M_0 \le M\cap k^{-1}Mk\cap kMk^{-1}.
\]
Therefore $M\cap k^{-1}Mk\cap kMk^{-1}$ is infinite, so Lemma~\ref{lem:almost_malnormal}(2) forces $\alpha|_M$ to be constant modulo $C_0(M)$, a contradiction.
\end{proof}

Our second obstruction does not require periodicity and eliminates almost malnormal defect hulls admitting elements whose conjugates almost all commute.

\begin{lemma}[Commuting-conjugates obstruction]\label{lem:commuting_conjugates}
Let $G$ be a group, let $\alpha\colon G\to\mathbb{C}$ be left-slowly oscillating, and let $M<G$ be an infinite defect hull for $\alpha$ that is almost malnormal in $G$. Suppose there exists $g\in G\setminus M$ such that $\{u\in M : [g,ugu^{-1}]\neq e\}$ is finite. Then $\alpha|_M$ is constant modulo $C_0(M)$.
\end{lemma}

\begin{proof}
Suppose for contradiction that $\alpha|_M$ is not constant modulo $C_0(M)$. Applying Lemma~\ref{lem:oscillation}(2) to $\alpha|_M$ and replacing the resulting terms by their inverses in $M$, we obtain distinct $c_1\neq c_2$ in $\mathbb{C}$ and sequences $(a_k')_{k\ge 1}$ and $(b_k')_{k\ge 1}$ in $M$ with pairwise distinct combined terms such that $\alpha((a_k')^{-1})\to c_1$ and $\alpha((b_k')^{-1})\to c_2$ as $k\to\infty$.

Set $E\coloneqq\{u\in M : [g,ugu^{-1}]\neq e\}$ and $Q\coloneqq M\cap g^{-1}Mg$; both $E$ and $Q$ are finite by hypothesis. We construct sequences $(a_n)_{n\ge 1}$ and $(b_n)_{n\ge 1}$ recursively as subsequences of $(a_k')$ and $(b_k')$. At step $n$, put $a_n\coloneqq a_n'$. Consider the set $W_n\coloneqq\{b\in M : ga_n^{-1}bg\in M\}$. If $W_n\neq\varnothing$, fix $b_0\in W_n$; then for every $b\in W_n$,
\[
  g a_n^{-1} b b_0^{-1} a_n g^{-1} = (g a_n^{-1} b g)(g a_n^{-1} b_0 g)^{-1}\in M,
\]
so $a_n^{-1} b b_0^{-1} a_n\in Q$ and hence $b\in a_n Q a_n^{-1} b_0$. Thus $|W_n|\le|Q|<\infty$ (vacuously when $W_n=\varnothing$), so $a_n E\cup W_n$ is finite. Since the terms of $(b_k')$ are pairwise distinct, we may choose $j_n>j_{n-1}$ (with $j_0\coloneqq 0$) such that $b_n\coloneqq b_{j_n}'\notin a_n E\cup W_n$.

Set $u_n\coloneqq a_n^{-1}b_n\in M$. By choice of $b_n$, we have $u_n\notin E$, so
\[
  (u_n g u_n^{-1})g = g(u_n g u_n^{-1}),
\]
and $gu_n g\notin M$. Define
\[
  x_n \coloneqq u_n^{-1} g a_n^{-1}
  \qquad\text{and}\qquad
  y_n \coloneqq u_n g b_n^{-1}.
\]
Because $g\notin M$ and $a_n,b_n\in M$, we have $ga_n^{-1},gb_n^{-1}\notin M$. Since $u_n^{\pm 1}\in M$ and $M$ is a defect hull, $\alpha(x_n)=\alpha(ga_n^{-1})$ and $\alpha(y_n)=\alpha(gb_n^{-1})$. As the elements $a_n^{-1}$ and $b_n^{-1}$ are pairwise distinct, left-slow oscillation for $g$ gives $\alpha(ga_n^{-1})-\alpha(a_n^{-1})\to 0$ and $\alpha(gb_n^{-1})-\alpha(b_n^{-1})\to 0$, whence
\[
  \alpha(x_n)\longrightarrow c_1
  \qquad\text{and}\qquad
  \alpha(y_n)\longrightarrow c_2
  \qquad (n\to\infty).
\]

On the other hand, using $u_n^{-1}a_n^{-1}=(b_n^{-1}a_n)a_n^{-1}=b_n^{-1}$ and commutativity of $g$ with $u_n g u_n^{-1}$, we compute
\[
  u_n(gx_n)
  = (u_n g u_n^{-1})g a_n^{-1}
  = g(u_n g u_n^{-1})a_n^{-1}
  = g u_n g b_n^{-1}
  = gy_n.
\]
Since $gu_n g\notin M$ and $b_n^{-1}\in M$, we have $gy_n\notin M$, and hence $gx_n=u_n^{-1}(gy_n)\notin M$. Because $u_n\in M$ and $M$ is a defect hull, it follows that
\[
  \alpha(gy_n) = \alpha(u_n(gx_n)) = \alpha(gx_n)
  \qquad\text{for all }n\ge 1.
\]

Finally, if $x_n=x_m$, then $u_n^{-1}g a_n^{-1}=u_m^{-1}g a_m^{-1}$, so $g(a_n^{-1}a_m)g^{-1}=u_n u_m^{-1}\in M$, giving $a_m\in a_n Q$. Since the $a_m$ are pairwise distinct, at most $|Q|$ indices $m$ satisfy $x_m=x_n$. Symmetrically, $y_n=y_m$ implies $b_m\in b_n Q$, so at most $|Q|$ indices $m$ satisfy $y_m=y_n$. Therefore both sequences $(x_n)$ and $(y_n)$ eventually leave every finite subset of $G$. Left-slow oscillation for $g$ then yields $\alpha(gx_n)-\alpha(x_n)\to 0$ and $\alpha(gy_n)-\alpha(y_n)\to 0$. Combining these with $\alpha(gx_n)=\alpha(gy_n)$ gives $\alpha(x_n)-\alpha(y_n)\to 0$, contradicting $c_1\neq c_2$.
\end{proof}

\begin{theorem}[Structure of periodic $\aleph_1$-counterexamples]\label{thm:aleph1_structure}
Let $G$ be a periodic group of cardinality $\aleph_1$ whose Higson corona $\partial G$ is not a singleton. Let $f\in C_h(G)$ be not constant modulo $C_0(G)$, and set $\alpha(x)\coloneqq f(x^{-1})$. Then $G$ admits a strictly increasing, continuous $\omega_1$-chain of countably infinite defect hulls
\[
  M_0 < M_1 < \cdots < M_\beta < \cdots \qquad (\beta<\omega_1)
\]
with $G=\bigcup_{\beta<\omega_1} M_\beta$ (where $M_\lambda=\bigcup_{\beta<\lambda} M_\beta$ for limit ordinals $\lambda<\omega_1$) such that for every $\beta<\omega_1$:
\begin{enumerate}
  \item $\alpha|_{M_\beta}$ is not constant modulo $C_0(M_\beta)$, and $\alpha$ is constant on every right $M_\beta$-coset other than $M_\beta$;
  \item $M_\beta$ is almost malnormal in $G$ (in particular, $N_G(M_\beta)=M_\beta$);
  \item for every $x\in G\setminus M_\beta$, the function $h\mapsto\alpha(x^{-1}hx)-\alpha(x)$ belongs to $C_0(M_\beta)$, and consequently $T_G(M_0,M_\beta)\subseteq M_\beta$;
  \item for every $g\in G\setminus M_\beta$, the set $\{u\in M_\beta : [g,ugu^{-1}]\neq e\}$ is infinite.
\end{enumerate}
\end{theorem}

\begin{proof}
Enumerate $G=\{g_\beta : \beta<\omega_1\}$. By Lemma~\ref{lem:oscillation}(1), choose a countably infinite subset $S\subseteq G$ witnessing that $\alpha$ is not constant modulo $C_0(G)$, and set $M_0\coloneqq\operatorname{cl}_\alpha(S)$. We construct countable defect hulls $(M_\beta)_{\beta<\omega_1}$ by transfinite recursion. Assuming $M_\beta$ is a countable defect hull, we have $G\setminus M_\beta\neq\varnothing$ because $|G|=\aleph_1$; let $y_\beta$ be the first element of $G\setminus M_\beta$ in the enumeration and set
\[
  M_{\beta+1} \coloneqq \operatorname{cl}_\alpha\bigl(M_\beta\cup\{g_\beta,y_\beta\}\bigr),
\]
which is a countable defect hull by Lemma~\ref{lem:defect_hull}. At a limit ordinal $\lambda<\omega_1$, set $M_\lambda\coloneqq\bigcup_{\beta<\lambda} M_\beta$; this is a countable union of countable defect hulls and is therefore a countable defect hull (if $h\in M_\lambda$, then $h\in M_\beta$ for some $\beta<\lambda$, so $D_h(\alpha)\subseteq M_\beta\subseteq M_\lambda$). Each $M_\beta$ is infinite since $S\subseteq M_0\subseteq M_\beta$, we have $M_\beta<M_{\beta+1}$ because $y_\beta\in M_{\beta+1}\setminus M_\beta$, and $G=\bigcup_{\beta<\omega_1} M_\beta$ because $g_\beta\in M_{\beta+1}$.

Since $S\subseteq M_0\subseteq M_\beta$, $\alpha|_{M_\beta}$ is not constant modulo $C_0(M_\beta)$, establishing (1). Properties (2), (3), and (4) now follow immediately from Lemmas~\ref{lem:almost_malnormal}(1), \ref{lem:self_normalizing}, and~\ref{lem:commuting_conjugates}, respectively.
\end{proof}

\begin{corollary}\label{cor:almost_malnormal}
Let $G$ be an uncountable periodic group that contains no proper countably infinite almost malnormal subgroup. Then $\partial G$ is a singleton.
\end{corollary}

\begin{proof}
If $|G|\ge\aleph_2$, this is Theorem~\ref{thm:aleph2}. If $|G|=\aleph_1$ and $|\partial G|>1$, then Theorem~\ref{thm:aleph1_structure}(2) produces a proper countably infinite almost malnormal subgroup $M_0<G$, a contradiction.
\end{proof}

The commuting-conjugates obstruction applies to the following transfinite wreath construction. Its separation from the earlier criteria is proved after the singleton assertion.

\begin{theorem}[Iterated restricted wreath products]\label{thm:wreath}
Let $F$ be a nontrivial finite group and let $M_0$ be a countably infinite locally finite group. Define a continuous $\omega_1$-chain of countable locally finite groups by setting
\[
  M_{\beta+1} \coloneqq F\wr M_\beta = \left(\bigoplus_{x\in M_\beta} F\right)\rtimes M_\beta
  \qquad (\beta<\omega_1)
\]
(where $M_\beta$ acts on the base group $B_\beta\coloneqq\bigoplus_{x\in M_\beta} F$ by left shifts) and $M_\lambda\coloneqq\bigcup_{\beta<\lambda} M_\beta$ at limit ordinals $\lambda<\omega_1$, and set $G\coloneqq\bigcup_{\beta<\omega_1} M_\beta$. Then $G$ is a locally finite group of cardinality $\aleph_1$ and its Higson corona $\partial G$ is a singleton. If $F$ is a finite $p$-group and $M_0$ is a locally finite $p$-group for the same prime $p$, then $G$ is also a locally finite $p$-group.
\end{theorem}

\begin{proof}
We first check the construction. At a successor stage, finitely many elements of $F\wr M_\beta$ involve finitely many lamp coordinates and generate a finite subgroup $P\le M_\beta$ in the acting group. The union of the $P$-orbits of those coordinates is a finite set $D$, so these elements lie in the finite group $F^D\rtimes P$. Thus local finiteness passes to successor stages and to unions at limit stages. Each $M_\beta$ is countably infinite, and the inclusion $M_\beta<M_{\beta+1}$ is proper because the base group is nontrivial. Hence $|G|=\aleph_1$. If $F$ and $M_0$ are $p$-groups, the same induction uses finite $p$-groups $F^D\rtimes P$, proving the last assertion.

Suppose for contradiction that $\partial G$ is not a singleton, and let $(L_\gamma)_{\gamma<\omega_1}$ be the continuous increasing $\omega_1$-chain of countable defect hulls produced by Theorem~\ref{thm:aleph1_structure}. Define strictly increasing sequences $(\alpha_n)_{n\ge 0}$ and $(\gamma_n)_{n\ge 0}$ in $\omega_1$ recursively: set $\alpha_0\coloneqq 0$, and given $\alpha_n$, choose $\gamma_n\ge\alpha_n$ with $M_{\alpha_n}\subseteq L_{\gamma_n}$ and then $\alpha_{n+1}>\gamma_n$ with $L_{\gamma_n}\subseteq M_{\alpha_{n+1}}$ (possible since each $M_\alpha$ and $L_\gamma$ is countable and both chains exhaust $G$). Then $\delta\coloneqq\sup_n\alpha_n=\sup_n\gamma_n$ is a limit ordinal in $\omega_1$, and continuity gives
\[
  M_\delta = \bigcup_{n=0}^{\infty} M_{\alpha_n} = \bigcup_{n=0}^{\infty} L_{\gamma_n} = L_\delta.
\]

Fix $e\neq a\in F$ and let $g\in B_\delta\setminus\{e\}\subseteq M_{\delta+1}\setminus M_\delta$ be the element supported at $\{e\}$ with $g(e)=a$. For every $u\in M_\delta\setminus\{e\}$, the conjugate $ugu^{-1}$ is supported at $\{u\}\neq\{e\}$, so $g$ and $ugu^{-1}$ have disjoint supports in $B_\delta$ and therefore commute. Thus $[g,ugu^{-1}]=e$ for every $u\in M_\delta=L_\delta$, contradicting Theorem~\ref{thm:aleph1_structure}(4). Hence $\partial G$ is a singleton.
\end{proof}

For example, taking $F=\mathbb Z/p\mathbb Z$ and $M_0=\bigoplus_{n\ge1}\mathbb Z/p\mathbb Z$ gives an explicit family of locally finite $p$-groups with singleton coronas.

\begin{proposition}[Separation from the preceding criteria]\label{prop:wreath_separation}
For the group $G$ in Theorem~\ref{thm:wreath}, every $M_\beta$ is almost malnormal in $G$, and $[N_1,N_2]\neq\{e\}$ for any two nontrivial normal subgroups $N_1,N_2\trianglelefteq G$. Moreover, $G$ satisfies none of the sufficient hypotheses in Theorems~\ref{thm:transporter}, \ref{thm:commuting_product}, \ref{thm:aleph2}, and~\ref{thm:nonperiodic_2gen}, or in Corollaries~\ref{cor:self_normalizing}, \ref{cor:hyperabelian}, and~\ref{cor:almost_malnormal}.
\end{proposition}

\begin{proof}
To show that every $M_\beta$ is almost malnormal in $G$, we first verify that $M_\beta$ is almost malnormal in $M_{\beta+1}=B_\beta\rtimes M_\beta$. Write elements of $M_{\beta+1}$ as $(b,m)$ with $b\in B_\beta$ and $m\in M_\beta$, with action $(u\cdot b)(z)=b(u^{-1}z)$ and multiplication $(b_1,u_1)(b_2,u_2)=(b_1(u_1\cdot b_2),u_1 u_2)$, identifying $M_\beta$ with $\{(e,u):u\in M_\beta\}$ and $B_\beta$ with $\{(b,e):b\in B_\beta\}$. Take $g=(b,m)\in M_{\beta+1}\setminus M_\beta$, so $b\neq e$ has nonempty finite support $\operatorname{supp}(b)\subseteq M_\beta$. Since $(e,m)\in M_\beta$, we have $M_\beta\cap g M_\beta g^{-1}=M_\beta\cap(b,e)M_\beta(b,e)^{-1}$. For $u\in M_\beta$,
\[
  (b,e)(e,u)(b^{-1},e) = \bigl(b(u\cdot b^{-1}),u\bigr),
\]
which lies in $M_\beta$ if and only if $u\cdot b=b$. In particular, $u\operatorname{supp}(b)=\operatorname{supp}(b)$. Fixing $x_0\in\operatorname{supp}(b)$, we obtain $u\in\operatorname{supp}(b)x_0^{-1}$, so $|M_\beta\cap gM_\beta g^{-1}|\le|\operatorname{supp}(b)|<\infty$. Replacing $g$ by $g^{-1}$ shows that $M_\beta$ is almost malnormal in $M_{\beta+1}$. Almost malnormality is transitive: if $A\le B\le C$ with $A$ almost malnormal in $B$ and $B$ almost malnormal in $C$, then for $c\in B\setminus A$ the intersection $A\cap c^{-1}Ac$ is finite directly, while for $c\in C\setminus B$ it is contained in the finite set $B\cap c^{-1}Bc$. Since almost malnormality is also preserved under increasing unions, transfinite induction shows that $M_\beta$ is almost malnormal in $M_\gamma$ for all $\beta<\gamma<\omega_1$, and hence in $G$.

Next, let $N_1,N_2\trianglelefteq G$ be nontrivial normal subgroups and choose $e\neq x_1\in N_1$ and $e\neq x_2\in N_2$. Pick $\beta<\omega_1$ with $x_1,x_2\in M_\beta$, and fix $e\neq a\in F$. Let $w_1\in B_\beta\le M_{\beta+1}$ take the value $a$ at $e\in M_\beta$ and $e$ elsewhere. Then $y\coloneqq[x_1,w_1]\in N_1\cap B_\beta$ has support $\{e,x_1\}$ and is therefore nontrivial. Now let $w_2\in B_{\beta+1}\le M_{\beta+2}$ take the value $a$ at $e\in M_{\beta+1}$ and $e$ elsewhere; then $v\coloneqq[y,w_2]\in N_1\cap B_{\beta+1}$ has support $\{e,y\}\subseteq B_\beta$ with $v(e)=a^{-1}$. Consider $[x_2,v]=(x_2 v x_2^{-1})v^{-1}\in[N_2,N_1]\cap B_{\beta+1}$. Because $x_2\in M_\beta\setminus\{e\}$ and $M_\beta\cap B_\beta=\{e\}$, we have $x_2\notin B_\beta$, so $x_2\notin\operatorname{supp}(v)$. Hence $[x_2,v](x_2)=(x_2\cdot v)(x_2)=v(e)=a^{-1}\neq e$, showing that $[N_1,N_2]\neq\{e\}$.

Since $G$ is periodic of cardinality $\aleph_1$, it does not satisfy the hypotheses of Theorems~\ref{thm:aleph2} and~\ref{thm:nonperiodic_2gen}. Since $M_0$ is almost malnormal in $G$ and countably infinite, $N_G(M_0)=M_0<G$, so $G$ does not satisfy the hypotheses of Corollaries~\ref{cor:self_normalizing} and~\ref{cor:almost_malnormal}. Moreover, let $H\subseteq G$ be any countable subset and choose $\beta<\omega_1$ with $H\subseteq M_\beta$. If $t\in T_G(M_0,H)$, then $t^{-1}M_0 t\subseteq H\subseteq M_\beta$, whence $M_0\subseteq M_\beta\cap tM_\beta t^{-1}$. Because $M_0$ is infinite and $M_\beta$ is almost malnormal in $G$, this forces $t\in M_\beta$. Thus $T_G(M_0,H)\subseteq M_\beta$ is countable for every countable $H\subseteq G$, so the transporter hypothesis fails.

If $G=AB$ with $[A,B]=\{e\}$, then $A$ and $B$ are both normal in $G$, so the preceding normal-subgroup property forces one of $A,B$ to be trivial; hence the commuting-product hypothesis fails. Taking $N_1=N_2$ in that property shows that $G$ has no nontrivial abelian normal subgroup. Consequently, $G$ is not virtually hyperabelian: otherwise $G$ would contain a hyperabelian normal subgroup $H_0\trianglelefteq G$ of finite index (taking the normal core of a finite-index hyperabelian subgroup); the first nontrivial term of an ascending normal abelian series of $H_0$ would be a nontrivial abelian normal subgroup $A\trianglelefteq H_0$ having finitely many $G$-conjugates $A_1,\ldots,A_r\trianglelefteq H_0$. Setting $P\coloneqq A_1\cdots A_r\trianglelefteq G$ and $I_J\coloneqq\bigcap_{j\in J} A_j$ for $J\subseteq\{1,\ldots,r\}$, induction on $k\ge 1$ shows that $\gamma_k(P)\subseteq\langle I_J : |J|=k\rangle$ (since $[I_J,A_m]=\{e\}$ for $m\in J$ and $[I_J,A_m]\subseteq I_{J\cup\{m\}}$ for $m\notin J$). Thus $\gamma_{r+1}(P)=\{e\}$, so $P$ is a nontrivial nilpotent normal subgroup of $G$ and $Z(P)$ is a nontrivial abelian normal subgroup of $G$, a contradiction. By \v{S}unkov's theorem \citep{Shunkov1972}, $G$ also contains no involution with finite centralizer, so neither hypothesis of Corollary~\ref{cor:hyperabelian} applies.
\end{proof}

Combining the results of Sections~\ref{sec:algebraic}--\ref{sec:aleph1} reduces the problem for all uncountable locally nilpotent groups to locally finite $p$-groups of cardinality $\aleph_1$.

\begin{corollary}[Reduction of locally nilpotent counterexamples]\label{cor:locally_nilpotent}
Let $G$ be an uncountable locally nilpotent group with $|\partial G|>1$. Then $G$ is locally finite, $|G|=\aleph_1$, and
\[
  G=A\times P,
\]
where $A$ is finite and $P$ is a locally finite $p$-group of cardinality $\aleph_1$ for some prime $p$, with $|\partial P|>1$. In particular, $P$ admits the continuous chain of countably infinite almost malnormal subgroups described in Theorem~\ref{thm:aleph1_structure}, including its transporter and commuting-conjugates restrictions.
\end{corollary}

\begin{proof}
Suppose that $|\partial G|>1$. By Theorem~\ref{thm:nonperiodic_2gen}, $G$ must be periodic, and by Theorem~\ref{thm:aleph2}, $|G|=\aleph_1$. A finitely generated periodic nilpotent group is finite, so $G$ is locally finite. Moreover, a periodic locally nilpotent group is the restricted direct product $G=\bigoplus_q G_q$ of its Sylow $q$-subgroups; see \citep[Theorems~5.2.18 and~12.1.1]{Robinson1996}. Since a restricted direct product of countably many countable groups is countable and $|G|=\aleph_1$, at least one Sylow subgroup $P\coloneqq G_p$ is uncountable (hence $|P|=\aleph_1$). Setting $A\coloneqq\bigoplus_{q\neq p} G_q$ gives $G=A\times P$. If $A$ were infinite, Theorem~\ref{thm:commuting_product} would force $\partial G$ to be a singleton; hence $A$ is finite. Since $[G:P]=|A|<\infty$, Proposition~\ref{prop:finite_index} implies $|\partial P|>1$, and applying Theorem~\ref{thm:aleph1_structure} to $P$ yields the stated $\omega_1$-chain.
\end{proof}

\section{Open questions}\label{sec:questions}

The results of Sections~\ref{sec:aleph2} and~\ref{sec:aleph1} localize the problem of singleton finitary Higson coronas for periodic groups entirely to cardinality $\aleph_1$, and for locally nilpotent groups entirely to locally finite $p$-groups of cardinality $\aleph_1$.

\begin{question}\label{q:lf_aleph1}
Does every locally finite group (in particular, every locally finite $p$-group) of cardinality $\aleph_1$ have a singleton finitary Higson corona? In view of Theorem~\ref{thm:aleph1_structure} and Corollary~\ref{cor:locally_nilpotent}, a prerequisite for a negative answer is the existence of a locally finite group $G$ of cardinality $\aleph_1$ admitting a continuous $\omega_1$-chain of countably infinite almost malnormal subgroups $(M_\beta)_{\beta<\omega_1}$ such that $\{u\in M_\beta : [g,ugu^{-1}]\neq e\}$ is infinite for every $\beta<\omega_1$ and every $g\in G\setminus M_\beta$. We do not know whether a locally finite group satisfying these necessary conditions exists. Their existence alone would not imply a nonsingleton corona, since no converse has been proved.
\end{question}

\begin{question}\label{q:burnside}
Is the cardinality threshold $|G|\ge\aleph_2$ in Theorem~\ref{thm:aleph2} sharp for general periodic groups? In particular, for sufficiently large odd exponents $n$, is the finitary Higson corona of the free Burnside group $B(\aleph_1,n)$ a singleton?
\end{question}

The connectedness of the corona in Question~\ref{q:burnside} is already known. An uncountable locally finite group is one-ended by Holt's theorem. A periodic group that is not locally finite is also one-ended: the alternatives in the splitting theorem would make it infinite virtually cyclic or give a nontrivial amalgamated product or HNN extension over a finite subgroup, each of which contains an element of infinite order; see \citep[Theorem~IV.6.10]{DicksDunwoody1989} and \citep[Theorems~1.2 and~1.4]{Cornulier}. Question~\ref{q:burnside} concerns whether this connected full corona is a singleton.

\section*{Statements and Declarations}
\noindent\textbf{Funding.} \FundingDisclosure\par\smallskip
\noindent\textbf{Competing interests.} \InterestDisclosure\par\smallskip
\noindent\textbf{Data availability.} This is a theoretical article. No datasets were generated or analyzed.\par\smallskip
\noindent\textbf{Use of AI tools.} The author used Gemini and ChatGPT to assist with reviewing proof details, organizing the exposition, editing the language, and preparing the \LaTeX{} source. The author is responsible for the mathematical arguments, references, and final text.\par

\end{document}